\newcommand{\del}{\partial}
\newcommand{\CR}{\hbox{{$\cal R$}}}
\newcommand{\CF}{\hbox{{$\cal T$}}}

\newcommand{\extd}{{\rm d}}

\newcommand{\N}{\mathbb{N}}
\newcommand{\Z}{\mathbb{Z}}

\newcommand{\R}{\mathbb{R}}
\newcommand{\C}{\mathbb{C}}

\newcommand{\isom}{{\cong}}
\newcommand{\eps}{{\epsilon}}
\newcommand{\tens}{\mathop{\otimes}}
\newcommand{\la}{{\triangleright}}\newcommand{\ra}{{\triangleleft}}

\newenvironment{proof}{\goodbreak\noindent{\bf Proof\quad}}{{\
 $\hbox{$\sqcup$}\llap{\hbox{$\sqcap$}}$}\bigskip }

\newcommand{\grav}{\mathsf{G}}
\newcommand{\pla}{{\C[x]\bicross_{\hbar,\grav}\C[p]}}
\newcommand{\copla}{{\C[\bar p]\bicross_{\frac{1}{\hbar},
\frac{\grav}{\hbar}}\C[\bar x]}}
\newcommand{\note}[1]{}

\newcommand{\rcross}{{\triangleright\!\!\!<}}

\newcommand{\cobicross}{{\triangleright\!\!\!\blacktriangleleft}}
\newcommand{\bicross}{{\blacktriangleright\!\!\!\triangleleft}}
\newcommand{\dcross}{{\bowtie}}

\newcommand{\rbiprod}{{\cdot\kern-.33em\triangleright\!\!\!<}}
\newcommand{\lbiprod}{{>\!\!\!\triangleleft\kern-.33em\cdot}}

\newcommand{\Ad}{{\rm Ad}}

\newcommand{\id}{{\rm id}}
\newcommand{\<}{\langle}
\renewcommand{\>}{\rangle}

\newcommand{\und}[1]{{\underline {#1}}}
\newcommand{\cb}{\mathfrak{b}} 
\newcommand{\cg}{\mathfrak{g}}

\renewcommand{\o}{{}_{\scriptscriptstyle(1)}}
\renewcommand{\t}{{}_{\scriptscriptstyle(2)}}
\renewcommand{\th}{{}_{\scriptscriptstyle(3)}}
\newcommand{\fo}{{}_{\scriptscriptstyle(4)}}
\newcommand{\fiv}{{}_{\scriptscriptstyle(5)}}
\newcommand{\six}{{}_{\scriptscriptstyle(6)}}
\newcommand{\sev}{{}_{\scriptscriptstyle(7)}}
\newcommand{\eig}{{}_{\scriptscriptstyle(8)}}
\newcommand{\lt}{{}_{{\scriptscriptstyle(\infty)}}}
\newcommand{\lo}{{}_{{\scriptscriptstyle(1)}}}
\newcommand{\rt}{{}_{{\scriptscriptstyle(1)}}}
\newcommand{\ro}{{}_{{\scriptscriptstyle(0)}}}
\newcommand{\uo}{{{}^{\scriptscriptstyle(1)}}}
\newcommand{\ut}{{{}^{\scriptscriptstyle(2)}}}
\newcommand{\uth}{{{}^{\scriptscriptstyle(3)}}}
\newcommand{\umo}{{{}^{\scriptscriptstyle-(1)}}}
\newcommand{\umt}{{{}^{\scriptscriptstyle-(2)}}}
\newcommand{\umth}{{{}^{\scriptscriptstyle-(3)}}}
\newcommand{\alt}{{}_{{\scriptscriptstyle<\infty>}}}
\newcommand{\alo}{{}_{{\scriptscriptstyle<1>}}}
\newcommand{\ao}{{}_{\scriptscriptstyle<1>}}
\newcommand{\at}{{}_{\scriptscriptstyle<2>}}
\newcommand{\ath}{{}_{\scriptscriptstyle<3>}}

\newcommand{\eqn}[2]{\begin{equation}#2\label{#1}\end{equation}}

\documentclass[11pt]{article}
\usepackage{amssymb,amsmath,amscd}

\newcommand{\bo}{{}_{\scriptscriptstyle(\underline 1)}}
\newcommand{\bt}{{}_{\scriptscriptstyle(\underline 2)}}
\newcommand{\bth}{{}_{\scriptscriptstyle(\underline 3)}}

\newcommand{\bfi}{{}_{\scriptscriptstyle(\underline 5)}}
\newcommand{\abo}{{}_{\scriptscriptstyle<\underline 1>}}
\newcommand{\abt}{{}_{\scriptscriptstyle<\underline 2>}}
\newcommand{\abth}{{}_{\scriptscriptstyle<\underline 3>}}

\newcommand{\ubp}{\ensuremath{U(\cb_+)}}
\DeclareMathOperator{\cop}{\Delta}
\newcommand{\tensor}{\otimes}
\DeclareMathOperator{\cou}{\epsilon}
\newcommand{\antip}{S}
\newcommand{\act}{\triangleright}
\newcommand{\ract}{\triangleleft}

\newcommand{\pdiff}[1]{\frac{\partial}{\partial #1}}
\newcommand{\pdiffr}[2]{\frac{\partial #1}{\partial #2}}
\newcommand{\cbp}{\ensuremath{C(B_+)}}
\newcommand{\diff}{\mathrm{d}}
\newcommand{\tchi}{{\chi'}}
\newcommand{\hchi}{{{\chi'}'}}
\DeclareMathOperator{\adl}{Ad_L}
\newcommand{\cfunct}{\mathcal{F}}
\newcommand{\cfc}{\cfunct_\chi}
\newcommand{\cgd}{\mathcal{G}_\chi}
\newcommand{\cgu}{\mathcal{G}^\chi}
\newcommand{\cnat}{c}
\newcommand{\cnatd}{{\cnat_\chi}}
\newcommand{\cnatu}{{\cnat^\chi}}
\newcommand{\tensh}{\tensor_H}
\newcommand{\tenshd}{\tensor_{H_\chi}}
\newcommand{\tenshu}{\tensor_{H^\chi}}
\newcommand{\catbibi}[1]{{}^{#1}_{#1}\mathcal{M}^{#1}_{#1}}
\newcommand{\catlc}[1]{{}^{#1}_{#1}{\dot\mathcal{M}}}
\newcommand{\catrc}[1]{\dot\mathcal{M}^{#1}_{#1}}
\renewcommand{\O}{\mathcal{O}}

\allowdisplaybreaks 

\newtheorem{lemma}{Lemma}[section]
\newtheorem{propos}[lemma]{Proposition}

\newtheorem{theorem}[lemma]{Theorem}

\newtheorem{corol}[lemma]{Corollary}

\begin{document}

{\ }\qquad  \hskip 4.3in DAMTP/98-118
\vspace{.2in}

\begin{center} {\LARGE TWISTING OF QUANTUM DIFFERENTIALS AND THE
PLANCK SCALE HOPF ALGEBRA}
\\ \baselineskip 13pt{\ }\\
{\ }\\ Shahn Majid\footnote{Royal Society University
Research Fellow and Fellow of Pembroke College, Cambridge}+
Robert Oeckl\\ {\ }\\
Department of Applied Mathematics \& Theoretical Physics\\
University of Cambridge, Cambridge CB3 9EW\\
\end{center}
\begin{center}
October, 1998
\end{center}

\begin{quote}
\noindent{\bf Abstract}
We show that the crossed modules and
bicovariant different calculi on two Hopf algebras related by a
cocycle twist are in 1-1 correspondence. In particular, for
quantum groups which are cocycle
deformation-quantisations of classical groups the calculi
are obtained as
deformation-quantisation of the classical ones. As an
application, we classify all bicovariant differential calculi on
the Planck scale Hopf algebra $\C[x]\bicross_{\hbar,\grav}\C[p]$.
This is a quantum group which has an $\hbar\to 0$ limit as the
functions on a classical but non-Abelian group and a $\grav\to 0$
limit as flat space quantum mechanics. We further study the
noncommutative differential geometry and Fourier theory for this Hopf
algebra as a toy model for Planck scale physics. The Fourier theory
implements a T-duality like self-duality. The noncommutative geometry
turns out to be singular when $\grav\to 0$ and is therefore not visible
in flat space quantum mechanics alone.


\end{quote}

\section{Introduction}

Recent years have seen considerable advances in the noncommutative
geometry related to quantum groups, including a more or less
complete theory of quantum bundles and connections with quantum
group structure, quantum homogeneous spaces, etc. Particularly
important for all these constructions is the differential calculus
or bimodule $\Omega^1$ of 1-forms on the `quantum space'. Recently,
the translation-bicovariant calculi on the quantum group itself
have been classified for the class of `factorisable' quantum groups
\cite{Ma:cla}. Using different methods, one also has a
classification of the calculi on bicrossproduct Hopf algebras
$\C(M)\bicross\C G$ associated to finite group factorisations
\cite{BegMa:dif}.

In this paper we present a third and novel `deformation theoretic'
approach to the construction of differential calculi which works
for quantum groups which are cocycle
deformation-quantisations \cite{Dri:hop} of classical groups. We
show that the calculi on the quantum group are in correspondence
with the calculi on the classical group and we provide the explicit
deformation-quantisation of the latter to obtain the former. This
includes the large class of quantum groups related to triangular
solutions of the Yang-Baxter equation \cite{Dri:hop}, not
covered by the previous approaches.

The paper consists of a general functorial result (in Section 3)
concerning how the construction of quantum differential calculi on
a quantum group responds under a certain `cotwisting' or `gauge
equivalence' \cite{Dri:qua}\cite{Ma:book} operation in the category
of quantum groups. The main idea of this cotwisting is to start
with an initial Hopf algebra $H$, e.g.\ a commutative one (the
coordinate ring of a classical group) and modify its product by a
cocycle $\chi$ to a new quantum group $H^\chi$. We show that the
exterior algebra of the quantum differential calculus likewise
twists as a super-Hopf algebra. This is the main result of
Section~2.

The main technical result is the construction of a nontrivial
monoidal functor relating the bicovariant bimodules over the
twisted and untwisted Hopf algebras. Moreover, bicovariant
bimodules over a Hopf algebra $H$ are equivalent to the
representations of the quantum double $D(H)$ or,
equivalently, to the crossed modules (Drinfeld-Radford-Yetter
modules) $\catlc{H}$ associated to $H$. Therefore a corollary, also
in Section~2, is the construction of a nontrivial monoidal functor
\[ \cfunct^\chi:\catlc{H}\to \catlc{H^\chi}\]
or equivalently (by Tannaka-Krein reconstruction) of an
isomorphism
\[ D(H^\chi)\isom D(H)_{\tilde\chi}\]
(for twisting of the coproduct by a certain other cocycle
$\tilde\chi$ in $D(H)$).

Finally, in Section~3, we add the consideration of the exterior
differential to complete the general theory. As a further
nontrivial corollary related to the particular crossed modules
needed for quantum differential calculi, we obtain an
identification (via the above functor) of the quantum adjoint
action under the twisted and untwisted Hopf algebras, extending a
result in \cite{GurMa:bra}.

The second half of the paper consists of the application of these
general results to a particular Hopf algebra
$\C[x]\bicross_{\hbar,\grav}\C[p]$ which was introduced 10 years
ago as a new `Hopf algebra approach' to Planck scale
physics \cite{Ma:pla}. It appears to be the first serious attempt to
develop gravitationally-modified quantum mechanics as
`noncommutative geometry' by quantum group methods, namely by
modifying the usual $[x,p]=\imath\hbar$ commutation relations in
such a way as to allow a coproduct. Thus, the Planck scale Hopf
algebra has two parameters $\hbar,\grav$ and the Hopf algebra
structure
\[
 [x,p]=\imath\hbar(1-e^{-\frac{x}{\grav}}),\quad \Delta x=x\tens 1+1\tens x,
\quad \Delta p=p\tens e^{-\frac{x}{\grav}}+1\tens p.
\]
In the limit
$\hbar\to 0$ one obtains the classical ring $C(B_+)$ of functions
on a group $B_+=\R\rcross_\grav
\R$ (a non-Abelian but classical group as the classical phase space)
and in another limit $\grav\to 0$ one obtains the usual flat space
Heisenberg algebra when restricted to $x>0$ (the coalgebra is
singular in this limit). Here $\hbar$ and $\grav$ are arbitrary
constants, but the latter plays a role similar to (a dimensionful
multiple of) the gravitational Newton coupling constant when one
compares free particle motion with the motion of a particle falling
into a black hole \cite{Ma:pla}\cite{Ma:book}. More precisely, $\grav$
plays the role of $\grav_{\rm Newton}M/c^2$ where $M$ is the mass of
the black hole and $c$ is the speed of light. Thus we can also
consider $\grav$ more physically as playing the role of the
gravitational mass or radius of curvature of the background
geometry.

A main feature of this Planck scale Hopf algebra is its
self-duality, i.e.\ the linear functionals on this quantum system
(containing the states) can be convolved and as such form an
isomorphic Hopf algebra which can be viewed as the observables of a
dual quantum system (the original algebra of observables containing
the states of this). It seems likely that more recent constructions
of $T$-duality in string theory can be viewed as generalisations to
string theory of such a duality. On the other hand, the
noncommutative differential geometry of the Planck scale Hopf
algebra could not be explored 10 years ago; we are able to do this
now. While we recover the known differential calculi on $B_+$
as $\hbar\to 0$, the differential calculi are equally valid for general 
$\hbar$ and thereby extend our geometrical notions to the quantum system. 
As well as the the classical and flat space quantum mechanical limits
there is a third limit where $\imath\kappa=\grav/\hbar$ is held fixed as
$\hbar,\grav\to \infty$, and the Planck scale Hopf algebra tends to
the enveloping algebra $U(\cb_-)$ ($\cb_-$ the opposite of the Lie algebra
$\cb_+$ of $B_+$). Such an algebra has been proposed as a noncommutative 
model of spacetime (here in 2 dimensions) covariant under the
$\kappa$-Poincar\'e quantum group \cite{MaRue:bic}; its noncommutative 
geometry can therefore be obtained as a special case.

It is known that the Planck scale Hopf algebra is a
cocycle twist of $U(\cb_+)$. In view of its self-duality this
implies that it is also a cotwist deformation quantisation of
$C(B_+)$, where $B_+$ is the group above. We begin this part of the
paper in Section~4 by obtaining the cotwisting cocycle $\chi_\hbar$
explicitly. In this way we have explicitly the
deformation-quantisation which undoes the $\hbar\to 0$ limit above
as a switching on of the cocycle $\chi_\hbar$. The cocycle
$\chi_\hbar$ has in fact a similar form to the operation in the
Moyal product \cite{BFFLS:def}, i.e.\ this is a variant of the Moyal
product or
$*$-product quantisation, the variation being that we do not start
from $\R^2$ but from a non-Abelian (i.e. in some sense `curved')
classical group manifold $B_+$. On the other hand, the first order
differential calculi on $C(B_+)$ have been completely classified
recently in \cite{Oe:cla} and hence our general construction in
Section~3 provides, in Section~4, a functorial
`deformation-quantisation' of these to calculi on the Planck scale
Hopf algebra.

We then compute the entire quantum exterior algebra in the case of
the quantisation of the standard classical 2-dimensional
differential calculus on $B_+$. We also provide some elements of
quantum Poisson theory in this case, including quantum-geometric
Hamilton equations of motion. Remarkably, the exterior algebra and
this quantum-geometrical picture is highly singular when $\grav\to
0$, i.e. is not visible in flat space quantum mechanics alone. As
in \cite{Ma:pla}, we conclude that the presence of even a small
amount of `gravity' makes quantum mechanics better behaved and
restores its geometry.

The remaining Section~5 completes the noncommutative picture with
formulae for the left and right invariant integrals, exponentials
and the Fourier transform on the Planck scale Hopf algebra. In fact
the quantum Fourier transform is a linear isomorphism
\[\CF: \C[x]\bicross_{\hbar,\grav}\C[p]\to
\C[\bar p]\cobicross_{\frac{1}{\hbar},\frac{\grav}{\hbar}}\C[\bar x],\]
which interchanges the roles of $x,p$ (a version of Born
reciprocity) and at the same time requires inversion of $\hbar$.
Thus we see that it explicitly implements the T-duality-like
feature of our toy model of Planck scale physics in \cite{Ma:pla}.
This self-duality isomorphism is singular when $\hbar\to 0$, i.e.\
only visible due to the presence of quantisation. When $\hbar\to 0$,
one obtains instead a non-Abelian Fourier isomorphism between two
completely different objects, namely
\[ \CF: C(B_+)\to U(\cb_+).\]
The right hand side here is also a version of $\kappa$-Minkowski space, with
$\imath\kappa=\grav$, so Fourier theory on this is recovered in the
classical limit.

\subsection*{Preliminaries}

The general results in Sections 2,3 work over a general field $k$.
The application results in Sections~4,5 are over $k=\C$. We use the
usual notations for a Hopf algebra $(H,\Delta,S,\eps)$ where $H$ is
a unital algebra, $\Delta:H\to H\tens H$ and $\eps:H\to k$ the
counital coalgebra and $S:H\to H$ the antipode or generalised
`inversion'. For convenience we require that $S$ is invertible. 
We use the Sweedler notation $\Delta h=h\o\tens h\t$
for $h\in H$ and similarly for a coaction $\beta:V\to H\tens V$ we
use the notation $\beta(v)=v\lo\tens v\lt$. Here coproducts and
coactions obey axioms like products and actions, but with
arrows reversed, see \cite{Ma:book} for an introduction.

Given a Hopf algebra $H$ and invertible $\chi\in H\tens H$ a cocycle
in the sense
\eqn{chi}{\chi_{23}(\id\tens\Delta\chi)=\chi_{12}(\Delta\tens\id)\chi,
\quad(\eps\tens\id)\chi=1,}
one has a new Hopf algebra $H_\chi$ with the same algebra structure
and counit as $H$ and the new coproduct and antipode
\eqn{twist}{ \Delta_\chi=\chi\Delta(\ )\chi^{-1},\quad S_\chi=US(\ )U^{-1},
\quad U=\chi\uo (S \chi\ut  ).}
This is the twist of $H$, see \cite{Dri:qua}\cite{GurMa:bra}.
Dually, given a Hopf algebra $H$ and $\chi:H\tens H\to k$ convolution
invertible, a cocycle in the sense
\eqn{cochi}{ \chi(g\o\tens f\o)\chi(h\tens g\t f\t)=\chi(h\o\tens g\o)
\chi(h\t g\t\tens f),
\quad\chi(1\tens h)=\eps(h)}
one has a new Hopf algebra $H^\chi$ with the same coalgebra and unit
as $H$ and the new product and antipode
\eqn{cotwist}{ h\bullet g=\chi(h\o\tens g\o)h\t g\t
 \chi^{-1}(h\th\tens g\th),
\quad S^\chi h= U(h\o)Sh\t U^{-1}(h\th),\quad U(h)=\chi(h\o\tens Sh\t).}
See \cite{Ma:book}.

Next, given a Hopf algebra (with invertible antipode), one has a
braided category $\catlc{H}$ of crossed modules 
\cite{Dri}\cite{Rad:str}\cite{Yet:rep}. Objects are vector
spaces $V$ which are both $H$-modules and $H$-comodules. The two
structures obey a compatibility condition. This and the
braiding $\Psi_{V,W}:V\tens W\to W\tens V$ are
\eqn{crossmod}{  h\o v\lo\tens h\t\la v\lt= (h\o\la v)\lo h\t\tens
(h\o\la v)\lt,\quad \Psi_{V,W}(v\tens w)= v\lo\la w\tens v\lt,}
where $\la$ denotes the action.
This is a slight reformulation (in a completely
standard manner) of the braided category  of modules over the
Drinfeld quantum double $D(H)$ (here $D(H)=H^{*\rm op}\dcross H$ 
generated when $H$ is finite-dimensional by $H$ and $H^{*\rm op}$
as sub-Hopf algebras and an $H^{*\rm op}$-module structure is 
equivalent to a $H$-comodule structure; the latter formulation then
avoids finite-dimensionality). The cross relations in $D(H)$ are
\eqn{double}{    h\phi= \phi\t h\t \<Sh\o,\phi\o\>\<h\th,\phi\th\>,
\quad \forall h\in H,
\quad\phi\in H^*.}
$H$ is a crossed module over itself with action by
left-multiplication and the adjoint coaction
\[ \beta(h)=h\o Sh\th\tens h\t.\]
Here, $\ker\eps$ is a sub-crossed module.

A differential calculus over any algebra $H$ means a specification of
an $H$-bimodule of differential 1-forms $\Omega^1$ and a map
$\extd:H\to\Omega^1$ obeying the Leibniz rule
\eqn{leib}{ \extd(hg)=(\extd h)g+h\extd g}
and surjective in the sense
$\Omega^1={\rm span}\{h\extd g|\ h,g\in H\}$. 
Such 1-forms can then be extended to an entire exterior algebra
with $\extd^2=0$, although not necessarily uniquely. When $H$ is a
Hopf algebra it is natural to require that $\Omega^1$ is covariant
under left and right coactions of $H$ \cite{Wor:dif}. Thus we require
$\Omega^1$ to be both a bimodule and a bicomodule, with the
coactions bimodule maps and mutually commuting, i.e.
$\Omega^1\in\catbibi{H}$, the category of {\em bicovariant bimodules}
over $H$. It is
known that bicovariant bimodules are equivalent to crossed modules
under the correspondence $\Omega^1=V\tens H$
for some $V\in \catlc{H}$. The bimodule structure on $\Omega^1$ is
\eqn{bimod}{ h\cdot(v\tens g)=h\o\la v\tens h\t g,\quad (v\tens g)\cdot h
=v\tens gh.}
The left and right coactions are likewise induced from $V$ and $H$
using the given coaction on $V$ and the coproduct of $H$. The
remaining axioms of a differential calculus then reduce to $V$ a
quotient of $\ker\eps\in \catlc{H}$. The corresponding $\extd$ is
\eqn{extd}{ \extd h=(\pi\tens\id)(h\o\tens h\t-1\tens h),}
where $\pi$ is the projection from $\ker\eps$ to $V$. The space $V$
is the space of right-invariant differential forms.

Moreover, given $\Omega^1$ there is its {\em
maximal prolongation} differential graded algebra
$\widetilde{\Omega^\cdot}=\oplus_n \Omega^n$ where
$\Omega^0=H$, $\Omega^1=\Omega^1$ and $\Omega^n$ is a
certain quotient of $\Omega^1\tens_H\cdots\tens_H\Omega^1$ ($n$-fold). The
product is given by $\tens_H$ or with degree 0 by the bimodule structure. The
differential structure is $\extd$ extended by $\extd^2=0$ and the
graded Leibniz rule. The quotienting in $\widetilde{\Omega^\cdot}$ is
the minimum required for the extension of $\extd$. There is also
a {\em Woronowicz prolongation} differential graded algebra
$\widehat{\Omega^\cdot}$ which is
likewise a quotient of the tensor algebra, this time using the
braiding induced from $\catlc{H}$ to
`skew-symmetrise'. It is a quotient of the maximal prolongation,
and for q-deformed examples is typically a
q-deformation of the classical exterior algebra.

Finally, dual to the vector space $V$ corresponding to a
bicovariant calculus on $H$ is a subspace $L\subseteq\ker\eps\subset
U$ of a Hopf algebra $U$ dual to $H$. We have a
self-contained notion of `space of right-invariant vector fields'
$L\subseteq\ker\eps\in\catlc{U}$, where $\ker\eps\subset U$ is a
crossed module by the coproduct and the quantum adjoint action. The
latter restricts to a map $L\tens L\to L$ with the result that $L$
is sometimes called a `quantum Lie algebra'.

Moreover, any $x\in L$ acts on $H$ by $\tilde{x}(h)=\<x,h\o\>h\t$, say.
One may view both $H$ (as above) and $L=V^*$ in the
braided category $\catlc{H}$, where $L$ is in
this category by
\eqn{Lcross}{ h\la x=\<Sh,x\o\>x\t,\quad \<\phi,x\lo\>x\lt
=(S^{-1}\phi\t)x\phi\o}
for $h\in H$ and $\phi\in U$ (the coregular action and
coadjoint coaction). Then one has the
braided-Leibniz rule
\eqn{braleib}{ \tilde{x}(hg)=(\tilde{x}(h))g+ h_i\, \widetilde{x^i}(g)}
where  $\Psi^{-1}(x\tens h)\equiv h_i\tens x^i$ (sum over the
index) is the inverse of the braiding between $L,H$, cf.
\cite{Ma:cla} in different conventions. Thus the elements of $L$
indeed act as braided derivations or `vector fields'. To work with
them we will employ the following convenient notation: if
$\{v,w,\cdots\}$ is a given basis of $V$ with dual basis
$\{v^*,w^*,\cdots\}$ say, we denote the braided vector field
$\widetilde{v^*}$ by $\del_{v}$, etc. One has $\extd
h=v\del_v(h)+w\del_w(h)+\cdots$ (sum over the basis).

\section{Twisting of bicovariant bimodules and crossed modules}

In this section we study the response of the
category of bicovariant bimodules and the category of crossed
modules of a Hopf algebra $H$ under twisting. We will later apply this (in
Section~3) to obtain the response under twisting of a quantum differential
calculus $\Omega^1$ and its exterior algebra.

\subsection{Bicovariant bimodules and tensor algebras}

Let $H$ be a Hopf algebra. We recall from the Preliminaries that a
bicovariant bimodule $\Omega$ means an $H$-bimodule and
$H$-bicomodule where the coactions are each bimodule maps. We
denote the coactions by
\[ \beta_L(\alpha)=\alpha\lo\tens\alpha\lt,\quad \beta_R(\alpha)
=\alpha\ro\tens\alpha\rt.\]
Since these coactions commute, we write
\[ \alpha\lo\tens\alpha\lt\ro\tens\alpha\lt\rt=\alpha\ro\lo
\tens\alpha\ro\lt\tens\alpha\ro
=\alpha\o\tens\alpha\bt\tens\alpha\th\]
when both are applied. Here the underline marking the component
living in $\Omega$ is an alternative (but unconventional) notation for
coactions which allows renumberings in the same manner as for
coproducts.

It is known (Brzezinski's theorem) that the Woronowicz exterior
algebra is a super-Hopf algebra \cite{Brz:rem}. This is also known
(the case given in detail in \cite{Brz:rem}) for a similar exterior
algebra where the quotient is generated in degree 2. We begin with an
analogous result for the tensor algebra on $\Omega$ in a
slightly more general form as a braided group. The proof is
analogous to these known results and hence we give it here only in
outline form. We also recall that the category of $\Z$-graded
spaces is braided with braiding
\[\Psi(\alpha\tens\beta)=q^{|\alpha||\beta|}\beta\tens\alpha\]
for any invertible $q$. It can be identified as the the category of
comodules under a quantum group $\Z_q$ \cite{Ma:book}.
We denote by $\Omega^{\tens_H n}$ the n-fold tensor product
$\Omega\tens_H\cdots\tens_H\Omega$.

\begin{propos}cf.\ \cite{Brz:rem}
\label{prop:tensor_algebra}
Let $\Omega$ be an $H$-bicovariant bimodule. The tensor algebra
$T_q\Omega=\oplus_n \Omega^{\tens_H n}$ is a braided group in the
category of $\Z_q$-comodules. The coproduct and antipode are
\[ \cop =\beta_L+\beta_R,\quad \antip \alpha=-(\antip \alpha\o)\cdot
\alpha\bt\cdot(\antip \alpha\th)\]
on degree 1 and extended to $T_q\Omega$ as a braided group.
\end{propos}
\begin{proof}
The proof is by induction. First note that $\Delta$ as stated is a
bimodule map since $\beta_L,\beta_R$ are. We extend it by
\[ \Delta(\alpha\tens_H\beta)
=q^{|\alpha\t||\beta\o|}\alpha\o\tens_H\beta\o\tens\alpha\t\tens_H\beta\t\]
which is well-defined since $\Delta$ on $\alpha,\beta$ is a
bimodule map. Moreover, for the same reason $\Delta$ remains a
bimodule map. Coassociativity on degree 1 follows from that of $H$
and the bicomodule properties of $\Omega$ (this step is the same
as in \cite{Brz:rem}), and likewise extends to all degrees by
induction. By construction, $\Delta$ is an algebra map with the
braided tensor product in the category of $\Z_q$-comodules. Hence we
have a bialgebra in the category of $\Z_q$-comodules.

Similarly, it is easy to see from $\beta_L,\beta_R$ bimodule maps
that $S(h\cdot\alpha)=(S\alpha)\cdot S h$ and $S(\alpha\cdot h)=(Sh)\cdot
S\alpha$ ($S$ a skew-bimodule map). We extend $S$ to higher products by
$S(\alpha\tens_H\beta)=q^{|\alpha||\beta|}(S\beta)\tens_H(S\alpha)$
which is therefore well defined and remains a skew-bimodule map.
That the antipode axiom is fulfilled then only has to be verified
on degree 1, and extends by induction to all degrees. This is
easily verified and is the same as in \cite{Brz:rem}.
\end{proof}

The structure of $T_q\Omega$ as a bicovariant bimodule can be
recovered from the bialgebra structure in the following way: The
module structure is given by the multiplication with one factor in
$H$. The comodule structure is given by the coproduct with
subsequent projection of the respective component to $H$ (the
degree zero part).

\begin{lemma}
\label{lem:ctwist_free}
Let $\chi\in H\tensor H$ be a counital 2-cocycle for $H$. Then $\chi$
extends to a graded counital 2-cocycle $\chi\in T_q\Omega\tensor
T_q\Omega$ by inclusion. The twisting of the coproduct yields
$(T_q\Omega)_\chi$ as a bicovariant bimodule over $H_\chi$ and a
$\Z_q$-braided group (a Hopf algebra in the category of
$\Z_q$-comodules). In particular, the degree 0 part of
$(T_q\Omega)_\chi$ is $H_\chi$.
\end{lemma}
\begin{proof}
The cocycle is embedded in degree 0 and the braiding with degree 0
in the category of $\Z_q$-comodules is the trivial one (independent
of $q$). Hence for this special type of cocycle the requirements to
make $T_q\Omega$ a $\Z_q$-braided group have the same form as for
bosonic Hopf algebras. This cocycle condition then reduces on
$\chi$ to the usual cocycle condition for the construction of
$H_\chi$. Also, the twisting preserves degree, and is thus a
functor in the category of $\Z_q$-comodules. As for $T_q\Omega$, we
recover the bicovariant bimodule structure of the twisted objects
by using the product and the coproduct with subsequent projection.
\end{proof}

The restriction of this Lemma to $\Omega\subset T_q\Omega$ provides
the desired twist of bicovariant bimodules. That this gives rise to
an isomorphism of categories is the following theorem.

\begin{theorem}
\label{thm:ctwist_bibimod}
Let $H$ be a Hopf algebra, $\chi\in H\tensor H$ a counital
2-cocycle. There is an isomorphism of braided categories
$\cgd:\catbibi{H}\to\catbibi{H_\chi}$. $\cgd$
leaves the actions unchanged and transforms the coactions according to
\[(\beta_L)_\chi=\chi\beta_L\chi^{-1},\qquad
 (\beta_R)_\chi=\chi\beta_R\chi^{-1}.\]
The monoidal structure is provided by the (identity) natural
transformation
\[
\cnatd:\cgd(V)\tenshd\cgd(W)\to\cgd(V\tensh W),\quad
v\tenshd w \mapsto v\tensh w\]
for all $V,W\in\catbibi{H}$.
\end{theorem}
\begin{proof}
In view of the proposition above, we write actions as
multiplications $av\equiv a\cdot v$, $a\bullet v\equiv a\cdot^\chi v$
(and similarly for the right actions). By restricting
Lemma~\ref{lem:ctwist_free} to degree 1 we already know that $\cgd$
maps bicovariant bimodules to bicovariant bimodules. To see that
$\cgd$ is indeed a functor we have to show that it maps morphisms
to morphisms. Let $f:V\to W$ be a morphism in $\catbibi{H}$. Since
$\cgd$ leaves the actions invariant, we just have to check that $f$
remains a comodule map. For the left coaction,
\[
(\beta_L)_\chi\circ f=\chi(\beta_L\circ f)\chi^{-1}
=\chi((\id\tensor f)\circ\beta_L)\chi^{-1}
=(\id\tensor f)\circ(\chi\beta_L\chi^{-1})
=(\id\tensor f)\circ(\beta_L)_\chi
.\]
Accordingly for the right coaction.

Next, we show that $\cgd$ is monoidal. The only non-trivial part
is to check that $c_\chi$ is a morphism in $\catbibi{H_\chi}$. 
We only do the proof for the
left action and left coaction; the right handed cases follow by
symmetry. For the left action this is the commutativity of the
diagram
\[\begin{CD}
H\tensor\cgd(V)\tenshd\cgd(W) @>\cdot >>
 \cgd(V)\tenshd\cgd(W)\\
@V\id\tensor\cnatd VV  @VV\cnatd V\\
H\tensor\cgd(V\tensh W) @>\cdot >> \cgd(V\tensh W)
\end{CD}\]
which is evident since
\[\cnatd(h(v\tenshd w))=\cnatd(h v\tenshd w)= h v\tensh w =h(v\tensh w)
 =h\cnatd(v\tenshd w).\]
For the left coaction the diagram is
\[\begin{CD}
\cgd(V)\tenshd\cgu(W) @>(\beta_L)_\chi >>
 H\tensor\cgd(V)\tenshd\cgd(W)\\
@V\cnatd VV  @VV\id\tensor\cnatd V\\
\cgd(V\tensh W) @>(\beta_L)_\chi>> H\tensor\cgd(V\tensh W)
\end{CD}\]
We write the twisted coaction as $(\beta_L)_\chi(v)=v\ao\tens
v\abt$ for clarity, and we denote a second copy of $\chi$ by
$\tchi$. Then,
\begin{align*}
(\beta_L)_\chi\circ\cnatd(v\tenshd w)&=(\beta_L)_\chi(v\tensh w)
 =\chi\beta_L(v\tensh w)\chi^{-1}\\
&=\chi\uo v\o w\o \chi\umo \tensor
  \chi\ut (v\bt \tensh w\bt )\chi\umt \\
&=\chi\uo v\o w\o \chi\umo \tensor
  \chi\ut v\bt \tensh w\bt \chi\umt \\
&=\chi\uo v\o \tchi\umo \tchi\uo
   w\o \chi\umo \tensor
  \chi\ut v\bt \tchi\umt \tchi\ut
  \tensh w\bt \chi\umt \\
&=\chi\uo v\o \tchi\umo \tchi\uo
   w\o \chi\umo \tensor
  \chi\ut v\bt \tchi\umt \tensh
  \tchi\ut  w\bt \chi\umt \\
&=v\ao w\ao\tensor
  v\abt\tensh w\abt\\
&=(\id\tensor\cnatd)(v\ao w\ao \tensor
  v\abt \tenshd w\abt )\\
&=(\id\tensor\cnatd)\circ(\beta_L)_\chi (v\tenshd w) .
\end{align*}
For the braiding we have to show that the following digram
commutes
\[\begin{CD}
\cgd(V)\tenshd\cgd(W) @>\Psi_\chi >>
 \cgd(W)\tenshd\chi\cgd(V)\\
@V\cnatd VV  @VV\cnatd V\\
\cgd(V\tensh W) @>\Psi >> \cgd(W\tensh V)
\end{CD}\]
Since $\cnatd$ is the identity transformation, this means that the
braiding in the untwisted and twisted category should be the same.
This is,
\begin{align*}
&\cnatd\circ\Psi_\chi(v\tenshd w)\\
&=
 v\ao w\abo \antip_\chi w\at
 \tensh \antip_\chi v\at v\abth w\ath\\
&= v\ao w\abo \antip_\chi(v\at w\at)
 \tensh v\abth w\ath\\
&= \chi\uo  v\o  w\bo\chi\umo
 \antip_\chi(\chi\ut  v\t  w\t  \chi\umt )
 \tensh \chi\uth v\bth w\th  \chi\umth\\
&= \chi\uo  v\o  w\bo \chi\umo  \tchi\uo
 \antip \tchi\ut  \antip(\chi\ut  v\t  w\t  \chi\umt )
 \antip\hchi\umo \hchi\umt
 \tensh \chi\uth v\bth w\th  \chi\umth\\
&= \chi\uo  v\o  w\bo \chi\umo  \tchi\uo
 \antip (\chi\umt \tchi\ut ) \antip(v\t  w\t )
 \tensh \antip(\hchi\umo \chi\ut )
 \hchi\umt  \chi\uth v\bth w\th  \chi\umth\\
&= \chi\uo  v\o  w\bo \chi\umo\o
 \antip (\chi\umo\t ) \antip(v\t  w\t )
 \tensh \antip(\chi\ut\o )
 \chi\ut\t  v\bth w\th  \chi\umt \\
&= v\o  w\bo \antip(v\t  w\t )
 \tensh v\bth w\th \\
&= v\o  w\bo \antip w\t
 \tensh \antip v\t  v\bth w\th \\
&= \Psi\circ\cnatd (v \tenshd w) .
\end{align*}
We wrote $\hchi$ for a third copy of $\chi$ and used the notation
\begin{align*}
\chi\uo \tensor\chi\ut \tensor\chi\uth
 &=(1\tensor\chi)\cdot(\id\tensor\cop)\chi ,\\
\chi\umo \tensor\chi\umt \tensor\chi\umth
 &=(\cop\tensor\id)\chi^{-1}\cdot(\chi^{-1}\tensor 1) .
\end{align*}
For the invertibility, it is clear from the formulae that twisting
by $\chi^{-1}$ after twisting by $\chi$ gives the original objects
and morphisms.
\end{proof}

Next we give a different set of results where the product rather
than the coproduct of $H$ is twisted. Note that the category of
bicovariant bimodules over $H$ does not depend symmetrically on the
product and coproduct of $H$ (for example, the tensor product in
the category is $\tens_H$). Hence the following theorem is not in
any simple way the dual of the one above. We start with the
corresponding lemma.

\begin{lemma}
\label{lem:ptwist_free}
Let $\chi:H\tensor H\to k$ be a unital 2-cocycle on $H$. Then
$\chi$ extends to  a graded unital 2-cocycle $\chi:T_q\Omega\tensor
T_q\Omega\to k$ by defining $\chi$ to be zero on elements of degree
$\ge 1$. Moreover, the twisted product yields $(T_q\Omega)^\chi$ as
a bicovariant bimodule over $H^\chi$ and a $\Z_q$-braided group. In
particular, the degree 0 part of $(T_q\Omega)^\chi$ is $H^\chi$.
\end{lemma}\begin{proof}
We observe that the counit vanishes on elements of degree $\ge 1$.
The conditions for $\chi$ to be a graded unital 2-cocycle on
$T_q\Omega$ thus reduce to the the conditions for $\chi$ to be a
unital 2-cocycle on $H$.
\end{proof}

\begin{theorem}
\label{thm:ptwist_bibimod}
Let $H$ be a Hopf algebra, $\chi:H\tensor H\to k$ a unital
2-cocycle. There is an isomorphism of braided categories
$\mathcal{G^\chi}:\catbibi{H}\to\catbibi{H^\chi}$.
$\mathcal{G}^\chi$ leaves the coactions unchanged and transforms
the actions according to
\begin{align*}
h\bullet v &=\chi(h\o \tensor v\o )\, h\t
 v\bt\, \chi^{-1}(h\th \tensor v\th ) ,\\
v\bullet h &=\chi(v\o \tensor h\o )\,
 v\bt h\t
 \,\chi^{-1}(v\th \tensor h\th) .
\end{align*}
The monoidal structure is provided by the natural transformation
\begin{align*}
\cnatu:\cgu(V)\tenshu\cgu(W)&\to\cgu(V\tensh W)\\
 v\tenshu w &\mapsto \chi(v\o \tensor w\o )\,
  v\bt\tensh w\bt\,
  \chi^{-1}(v\th \tensor w\th )
\end{align*}
for all $V,W\in\catbibi{H}$.
\end{theorem}
\begin{proof}
By restricting Lemma \ref{lem:ptwist_free} to degree 1 we already
know that $\cgu$ maps bicovariant bimodules to bicovariant
bimodules. To see that $\cgu$ is indeed a functor we have to show
that it maps morphisms to morphisms. Let $f:V\to W$ be a morphism
in $\catbibi{H}$. Since $\cgu$ leaves the coactions invariant, we
just have to check that $f$ remains a module map. For all $h\in H$,
$v\in V$,
\begin{align*}
f(h\bullet v) &=\chi(h\o \tensor v\o )\, f(h\t  v\bt)
  \,\chi(h\th \tensor v\th )\\
 &=\chi(h\o \tensor v\o )\, h\t  f(v\bt)
  \,\chi(h\th \tensor v\th )\\
 &=\chi(h\o \tensor (f(v))\o )\, h\t  (f(v))\bt
  \,\chi(h\th \tensor (f(v))\th )\\
 &=h\bullet f(v).
\end{align*}
Similarly for the right action.

Next, we show that $\cgu$ is monoidal. The associativity property of $c^\chi$
and invertibility follow from $\chi$ an invertible cocycle. We verify
that $c^\chi$ is a morphism in $\catbibi{H^\chi}$. As before, we 
only do the proof for the left action and left coaction; the right 
handed versions follow by symmetry. For the left coaction this is
the commutativity of the diagram
\[\begin{CD}
\cgu(V)\tenshu\cgu(W) @>(\beta_L)^\chi >>
 H\tensor\cgu(V)\tenshu\cgu(W)\\
@V\cnatu VV  @VV\id\tensor\cnatu V\\
\cgu(V\tensh W) @>(\beta_L)^\chi>> H\tensor\cgu(V\tensh W)
\end{CD}\]
Explicitly,
\begin{align*}
(\beta_L)^\chi\circ\cnatu(v\tenshu w)
 &=(\beta_L)^\chi\left(\chi(v\o \tensor w\o )\,
    v\bt\tensh w\bt
   \,\chi^{-1}(v\th \tensor w\th )\right)\\
 &= \chi(v\o \tensor w\o )\,
   v\t  w\t \tensor
   v\bth\tensh w\bth
   \,\chi^{-1}(v\fo \tensor w\fo )\\
 &= \chi(v\o \tensor w\o )\,
   v\t  w\t \,\chi^{-1}(v\th \tensor w\th )\\
  &\quad \tensor
   \chi(v\fo \tensor w\fo )\,
   v\bfi\tensh w\bfi
   \,\chi^{-1}(v\six \tensor w\six )\\
 &=(\id\tensor\cnatu)(v\o \bullet w\o \tensor
   v\bt\tenshu w\bt)\\
 &=(\id\tensor\cnatu)\circ(\beta_L)^\chi(v\tenshu w) .
\end{align*}
For the left action we require the commutativity of the diagram
\[\begin{CD}
H\tensor\cgu(V)\tenshu\cgu(W) @>\bullet >>
 \cgu(V)\tenshu\cgu(W)\\
@V\id\tensor\cnatu VV  @VV\cnatu V\\
H\tensor\cgu(V\tensh W) @>\bullet >> \cgu(V\tensh W)
\end{CD}\]
This is,
\begin{align*}
\cnatu(h\bullet(v\tenshu w)) &=\cnatu(h\bullet v\tenshu w)\\
 &=\chi(h\o \tensor v\o )\,
   \cnatu(h\t  v\bt\tenshu w)\,\chi^{-1}(h\th
   \tensor v\th )\\
 &=\chi(h\o \tensor v\o )\,\chi(h\t  v\t \tensor w\o )\,
   h\th  v\bth\tensh w\bt\\
 &\quad  \chi^{-1}(h\fo  v\fo \tensor w\th )
   \,\chi^{-1}(h\fiv  \tensor v\fiv )\\
 &=\chi(h\o \tensor v\o )\,
   \chi^{-1}(h\t\o\tensor v\t\o)\\
 &\quad
   \chi(v\t\t\tensor w\o\o)\,
   \chi(h\t\t\tensor v\t\th w\o\t)\,
   h\th  v\bth\tensh w\bt\\
 &\quad  \chi^{-1}(h\fo  v\fo \tensor w\th )
   \,\chi^{-1}(h\fiv  \tensor v\fiv )\\
 &=\chi(v\o \tensor w\o )\,
   \chi(h\o \tensor v\t  w\t )\,
   h\t  v\bth\tensh w\bth\\
 &\quad  \chi^{-1}(h\th  v\fo \tensor w\fo )
   \,\chi^{-1}(h\fo  \tensor v\fiv )\\
 &=\chi(v\o \tensor w\o )\,
   \chi(h\o \tensor v\t  w\t )\,
   h\t  v\bth\tensh w\bth\\
 &\quad  \chi^{-1}(h\th\o\tensor v\fo\o w\fo\o)\,
   \chi^{-1}(v\fo\t\tensor w\fo\t)\\
 &\quad  \chi(h\th\t\tensor v\fo\th)\,
   \chi^{-1}(h\fo  \tensor v\fiv )\\
 &=\chi(v\o \tensor w\o )\,
   \chi(h\o \tensor v\t  w\t )\,
   h\t  v\bth\tensh w\bth\\
 &\quad  \chi^{-1}(h\th \tensor v\fo  w\fo )\,
   \chi^{-1}(v\fiv \tensor w\fiv )\\
 &=\chi(v\o \tensor w\o )\,
   \chi\left(h\o \tensor (v\bt\tensh w\bt)\o \right)\,
   h\t  (v\bt\tensh w\bt)\bt\\
 &\quad  \chi^{-1}\left(h\th \tensor (v\bt
   \tensh w\bt)\th \right)\,
   \chi^{-1}(v\th \tensor w\th )\\
 &=h\bullet\left(\chi(v\o \tensor w\o )\,
   v\bt\tensh w\bt\,
   \chi^{-1}(v\th \tensor w\th )\right)\\
 &=h\bullet\cnatu(v\tenshu w) .
\end{align*}

Next we show that $\cgu$ preserves the braiding. This is the
commutativity of the diagram
\[\begin{CD}
\cgu(V)\tenshu\cgu(W) @>\Psi^\chi >>
 \cgu(W)\tenshu\chi\cgu(V)\\
@V\cnatu VV  @VV\cnatu V\\
\cgu(V\tensh W) @>\Psi >> \cgu(W\tensh V)
\end{CD}\]
The braiding is $\Psi(h v\tensh w g)=h w\tens_H v g$
with $h,g\in H, v\in V$ left-invariant and $w\in W$
right-invariant (similarly over $H^\chi$). Thus we have,
\begin{align*}
\cnatu\circ\Psi^\chi(h\bullet v\tenshu w\bullet g)
 &=\cnatu(h\bullet w\tenshu v\bullet g)\\
 &=\cnatu\left(\chi(h\o \tensor w\o )\, h\t  w\bt
  \tenshu  v\bo g\o  \,
  \chi^{-1}( v\t \tensor g\t )\right)\\
 &=\chi(h\o \tensor w\o )\,
   \chi(h\t  w\t  \tensor g\o )\,h\th  w\bth
  \tensh  v\bo g\t  \\
 &\quad \chi^{-1}(h\fo \tensor v\t  g\th )\,
  \chi^{-1}( v\th \tensor g\fo )\\
 &=\Psi(\chi(h\o \tensor w\o )\,
   \chi(h\t  w\t  \tensor g\o )\,h\th  v\bo
  \tensh  w\bth g\t  \\
 &\quad \chi^{-1}(h\fo \tensor v\t  g\th )\,
  \chi^{-1}( v\th \tensor g\fo ))\\
 &=\Psi(\chi( w\o \tensor g\o )\,
   \chi(h\o \tensor w\t  g\t )\,h\t  v\bo
  \tensh  w\bth g\th  \\
 &\quad \chi^{-1}(h\th \tensor v\t  g\fo )\,
  \chi^{-1}( v\th \tensor g\fiv ))\\
 &=\Psi(\chi( w\o \tensor g\o )\,
   \chi(h\o \tensor w\t  g\t )\,h\t  v\bo
  \tensh  w\bth g\th  \\
 &\quad \chi^{-1}(h\th  v\t \tensor g\fo )\,
  \chi^{-1}(h\fo \tensor v\th ))\\
 &=\Psi\circ\cnatu(\chi( w\o \tensor g\o )\,
   h\o  v\bo \tenshu  w\bt g\t \,
  \chi^{-1}(h\t \tensor v\t ))\\
 &=\Psi\circ\cnatu(h\bullet v \tenshu w\bullet g) .
\end{align*}

Only the invertibility remains to be shown. The inverse
operation to twisting by $\chi$ is twisting by $\chi^{-1}$. Since both,
the twisting of the actions and the natural transformation $\cnatu$ look
formally exactly like the twisting of the product in $H$, we see that
applying $\chi^{-1}$ after $\chi$ will give the original object
in the same way as for the product in $H$.
\end{proof}

From these results we obtain in particular the twisting of the
bicovariant bimodule $\Omega$ itself to $\Omega_\chi$ and
$\Omega^\chi$ in  the two cases. Finally, we are able to identify
the braided groups in the above lemmas,

\begin{corol}
\label{cor:com_twist_star}
$(T_q\Omega)_\chi=T_q(\Omega_\chi)$ and $(T_q\Omega)^\chi\isom
T_q(\Omega^\chi)$ in the settings of Theorem~\ref{thm:ctwist_bibimod}
and~\ref{thm:ptwist_bibimod}
respectively.
\end{corol}
\begin{proof}
Given a bicovariant bimodule $\Omega$ over a Hopf algebra $H$ we
can identify $(T_q\Omega)_\chi$ with $(T_q\Omega_\chi)$ using the
natural transformation $\cnatd$ in Theorem~\ref{thm:ctwist_bibimod}
extended to multiple
tensor products. Similarly for $()^\chi$ using
Theorem~\ref{thm:ptwist_bibimod}. Since
$\cnat$ preserves degree, this identification is graded (a morphism
in the category of $\Z_q$-comodules).
\end{proof}

\subsection{Twisting of crossed modules}

In this section we restrict the above results to the
right-invariant part of the bicovariant bimodule $\Omega$. This is
a crossed module $V$ and $\Omega=V\tens H$ as explained in the
preliminaries. In this way we obtain, as corollaries of the
preceeding subsection, results about the response of crossed modules
under twisting. This is relevant to the our treatment of differential
calculi but it is also of independent interest in several other
algebraic settings where crossed modules play an important role.
For this reason some more direct proofs are provided in the
appendix.

\begin{theorem}
\label{thm:ptwist_lccat}
Let $H$ be a Hopf algebra, $\chi:H\tensor H\to k$ a unital
2-cocycle. There is an isomorphism of braided categories
$\cfunct^\chi:\ \catlc{H}
\to \catlc{H^\chi}$ given by the
identity on the underlying vector spaces and coactions, and
 transforming the action $\la$ to
\[h\act^\chi v=\chi(h\o\tensor v\lo) (h\t\act
v\lt)\lt \chi^{-1}((h\t\act v\lt)\lo
\tensor h\th) .
\]
The monoidal structure is given by the natural transformation
\[
\cnatu:\cfunct^\chi(V)\tensor\cfunct^\chi(W) \to
                \cfunct^\chi(V\tensor W),\quad
v\tensor w \mapsto \chi(v\lo\tensor w\lo) v\lt\tensor w\lt.
\]
\end{theorem}
\begin{proof}
We deduce this from Theorem~\ref{thm:ptwist_bibimod} using the
equivalence of categories $\catbibi{H}\cong\catlc{H}$. As explained
in the Preliminaries, a bicovariant bimodule has the canonical form
$\Omega=V\tensor H$ with $V$ a crossed module. Conversely, $V$ may
be recovered as the space of right-invariant forms of $\Omega$ with
the action $\la$ recovered from the bimodule structure on $\Omega$
(denoted $\cdot$) by
\[
h\act v=h\o \cdot v\cdot\antip h\t .
\]
Since the twisting in
$\catbibi{H}$ preserves the coactions, it preserves the
decomposition $\Omega=V\tensor H$ and thus induces a twisting
$\catlc{H}\to\catlc{H^\chi}$ by restriction to the right-invariant
forms $V$. The coaction of the crossed module remains
unchanged, while the twisted action in $\catlc{H^\chi}$ becomes (we
denote the twisted actions by $\act^\chi$ and $\bullet$
respectively):
\begin{align*}
h\act^\chi v &= h\o \bullet v\bullet\antip^\chi h\t \\
 &= h\o \bullet v\bullet\antip h\th  U(h\t )\, U^{-1}(h\fo )\\
 &= \chi(h\o \tensor v\o )\,
   h\t \cdot v\bt\bullet\antip h\fo \,
   U(h\th ) U^{-1}(h\fiv )\\
 &= \chi(h\o \tensor v\o )\,
    \chi(h\t  v\t \tensor\antip h_{(8)})\,
   h\th \cdot v\bth\cdot\antip h_{(7)}
   \chi^{-1}(h\fo  \tensor\antip h\six )\,
   U(h\fiv ) U^{-1}(h_{(9)})\\
 &= \chi(h\o \tensor v\o )\,
    \chi(h\t  v\t \tensor\antip h\fiv )\,
   h\th \cdot v\bth\cdot\antip h\fo \,
   U^{-1}(h\six )\\
 &= \chi(h\o \tensor v\o )\,
    \chi(h\t  v\t \tensor\antip h\fo )\,
   h\th \act v\bth\, U^{-1}(h\fiv )\\
 &= \chi(h\o \tensor v\o )\,
    \chi((h\t \act v\bt)\o  h\th \tensor\antip h\fo )\,
   (h\t \act v\bt)\bt\, U^{-1}(h\fiv )\\
 &= \chi(h\o \tensor v\o )\,
   (h\t \act v\bt)\bt\,
    \chi^{-1}((h\t \act v\bt)\o  \tensor h\th).
\end{align*}
We used the identity $\chi(a h\o \tensor \antip h\t )\,U^{-1}(h\th
)=\chi^{-1}(a\tensor h)$
and convert to a more conventional comodule notation as stated.
Finally, for the monoidal structure,
restricting $\cnatu$ given in Theorem~\ref{thm:ptwist_bibimod} to
the right-invariant component leads to the stated form. The
fact that $\cfunct^\chi$ preserves the braiding just follows
from the fact that the braiding in $\catlc{H}$ is induced by the
braiding in $\catbibi{H}$. The isomorphism property follows from
the categorial equivalence $\catbibi{H}\cong\catlc{H}$. \end{proof}

We proceed to give the dual version of this theorem.
Note that it can not be derived from Theorem~\ref{thm:ctwist_bibimod}
in a way analogous to the above proof.
The reason is essentially that
Theorems \ref{thm:ptwist_bibimod} and \ref{thm:ctwist_bibimod} are not
strictly dual to each other. This is because the tensor product
$\tensor_H$ is not self-dual in our sense. Dualisation instead
converts this tensor product to the corresponding cotensor product.
Since we do not want to concern ourselves with that here, we give a
proof by dualisation of Theorem~\ref{thm:ptwist_lccat}. In the crossed
module setting this presents no
further problem since the tensor product in $\catlc{H}$ is the usual
one and thus self-dual in our sense.

\begin{theorem}
\label{thm:ctwist_lccat}
Let $H$ be a Hopf algebra, $\chi\in H\tensor H$ a counital 2-cocycle
for $H$. There is an isomorphism of braided categories
${\cfunct}_\chi:\ \catlc{H}
\to \catlc{H_\chi}$ given by the
identity on the underlying vector spaces and actions,  and transforming
the coaction $\beta(v)=v\lo\tens v\lt$ to
\[\beta_\chi(v)=\chi\uo (\chi\umo \act v)\lo \chi\umt
 \tensor\chi\ut \act (\chi\umo \act v)\lt. \]
The monoidal structure is given by the natural transformation
\[
\cnatd:{\cfunct}_\chi(V)\tensor{\cfunct}_\chi(W) \to
                {\cfunct}_\chi(V\tensor W),\quad v\tensor w \mapsto
                \chi\umo \act v\tensor \chi\umt  \act w.
\]
\end{theorem}
\begin{proof}
This Theorem is strictly dual to Theorem~\ref{thm:ptwist_lccat} and
therefore equivalent:
We can write the action of Theorem~\ref{thm:ptwist_lccat} as
\begin{align*}
\act^\chi& = (\chi^{-1}\tensor\id)\circ(\id\tensor\tau)
 \circ(\beta\tensor\id)\circ(\act\tensor\id)\circ
 (\id\tensor\tau)\\
 &\quad\circ
 (\cop\tensor\id)\circ(\chi\tensor\id\tensor\id)
 \circ(\id\tensor\tau\tensor\id)\circ(\cop\tensor\beta)
\end{align*}
with $\chi:H\tensor H\to k$ and $\tau$ the flip map. Dualising
means reversing the order of the composition, exchanging product
with coproduct, action with coaction and switching to $\chi: k\to
H\tensor H$:
\begin{align*}
\beta_\chi &=(\cdot\tensor\act)\circ(\id\tensor\tau\tensor\id)\circ
 (\chi\tensor\id\tensor\id)\circ(\cdot\tensor\id)\\
 &\quad\circ
 (\id\tensor\tau)\circ(\beta\tensor\id)\circ
 (\act\tensor\id)\circ(\id\tensor\tau)\circ(\chi^{-1}\tensor\id).
\end{align*}
This is just the coaction stated. For the monoidal structure, we have to
take into account that $\cnat$ flips its direction under
dualisation. So the dual of
\[\cnatu=(\chi\tensor\id\tensor\id)
 \circ(\id\tensor\tau\tensor\id)\circ(\beta\tensor\beta)\]
in Theorem~\ref{thm:ptwist_lccat} is
\[(\cnatd)^{-1}=(\act\tensor\act)\circ(\id\tensor\tau\tensor\id)
 \circ(\chi\tensor\id\tensor\id).\]
Inverting replaces $\chi^{-1}$ with $\chi$, leading to the formula
stated.
\end{proof}

\begin{corol}
\label{cor:twist_double}
Let $H$ be finite-dimensional and $D(H)$ its quantum
double. Then there is a Hopf algebra isomorphism $\theta:D(H^\chi)\isom
D(H)_{\tilde\chi}$ where
$\tilde\chi=\chi^{-1}$ viewed in $D(H)\tens D(H)$, and the twisting is
that of the coproduct of $D(H)$. Here,
\[ \theta(\phi\tens h)=\chi\umt\o\tens\chi'\umo\phi\tens h\t
\<h\o,US\chi\umo\>\<h\th,\chi\umt\t\chi'\umt\>,\quad\forall
\phi\in H^*,\ h\in H.\]
\end{corol}
\begin{proof} This follows by Tannaka-Krein arguments \cite{Ma:book}.
Thus, there is a uniquely determined algebra isomorphism $\theta$
such that $\cfunct^\chi$ is pull-back along $\theta$, and
this is a bialgebra map up to a conjugation corresponding to the
nontrivial isomorphism of tensor products of objects in the image
of $\cfunct^\chi$ in Theorem~\ref{thm:ptwist_lccat}. We build
$D(H)$ explicitly on $H^*\tens H$ as explained in the
Preliminaries. Then,
\begin{align*}
\theta(\phi\tens h)\la v &=(\phi\tens h)\la^\chi v\\
&=\phi((h\t\la v\lt)\lt\lo) (h\t\la v\lt)\lt\lt
\chi(h\o\tens v\lo)\chi^{-1}((h\t\la v\lt)\lo\tens h\th))\\
&=\phi((h\t\la v\lt)\lo\t) (h\t\la v\lt)\lt\chi(h\o\tens v\lo)
\chi^{-1}((h\t\la v\lt)\o\tens h\th))\\
&=\phi\cdot\chi\umo\la(h\t\la v\lt)\<\chi\umt,h\th\>\chi(h\o\tens
v\lo)\\ &=\phi\cdot\chi\umo\cdot h\t\cdot\chi\ut\la
v\<\chi\umt,h\th\>
\<\chi\uo,h\o\>\\
&=\phi\cdot\chi\umo\cdot\chi\ut\t\cdot h\t\t\la v\<\chi\umt,h\th\>
\<\chi\uo,h\o\>\<S\chi\ut\o,h\t\o\>\<\chi\ut\th,h\t\th\>\\
&=\phi\cdot\chi\umo\cdot\chi\ut\t\cdot h\t\la v \<h\o,\chi\uo
S\chi\ut\o\>
\<h\th,\chi\ut\th\chi\umt\>
\end{align*}
using the definition of $\la^\chi$, the definition of the action of $D(H)$ on
$v\in V$ (namely $\phi$ acts by evaluation against the coaction),
then using the
cross relations of $D(H)$ and the duality pairing axioms. Here $\cdot$
denotes the
product in $D(H)$. From this and from
the cocycle identity
\[ \chi\uo S\chi\ut\o\tens \chi\ut\t=\chi\uo (S\chi\ut)
S\chi\umo\tens \chi\umt\]
(which follows from $\id \tens S$ and the product applied to
$(\id\tens\id)\chi
=\chi^{-1}_{23}\chi_{12}(\Delta\tens\id)\chi$), we find that
\[ \theta(\phi\tens h)=\chi\ut\t\chi\umo\phi\tens h\t\la v
\<h\o,\chi\uo S\chi\ut\o\>
\<h\th,\chi\ut\th\chi\umt\>\]
(where the products are in $H^*$ and $H$) has the form stated.
Next, we note that if two Hopf algebras have their coproducts
related by twisting by $\chi$ then the induced monoidal functor has
the same form as for $\cfunct^\chi$ on tensor products
\cite{Ma:book}. From this, we conclude that $\theta$ is an
isomorphism to $D(H)_{\tilde\chi}$ where $\tilde\chi
=\chi^{-1}\in H^*\tens H^*$ is viewed in $D(H)\tens D(H)$. Since $H^{*\rm op}
\subseteq D(H)$
as a subalgebra, conjugation by $\chi$ using the product of $H^*$ is actually
conjugation by $\chi^{-1}$ in $D(H)$. That $\theta$ is then an isomorphism of
bialgebras can also be checked explicitly.
\end{proof}

Clearly, the dual of the double $D(H)^*$ changes to
$D(H)^*{}^{\tilde\chi}$ similarly. For completeness, let us mention
also the dual version of Corollary~\ref{cor:twist_double}. Let
$\chi\in H\tens H$ be a
cocycle and $H_\chi$ the Hopf algebra with twisted coproduct.
Since the double $D(H)$ involves both $H$ and $H^{*\rm op}$
equally, $D(H_\chi)$ is likewise a twisting.

\section{Twisting of differential calculi and exterior algebras}

In this section we apply the technical results of Section~2 to obtain
our main result, which is a twisting theory for first order and exterior
differential calculi. We fix a Hopf algebra $H$ and
a first order differential calculus $(\Omega^1,\extd)$.
We begin with a more explicit discussion of the maximal prolongation exterior
algebra that we have found elsewhere, and the analogue of
Brzezinski's theorem \cite{Brz:rem} for it. We will then study how it,
and the more well-known Woronowicz exterior algebra, respond under twisting.

\begin{propos}
\label{prop:max_prolong}
The maximal prolongation $\widetilde{\Omega^\cdot}=\oplus_n\Omega^n$ is a
quotient of $T_{-1}\Omega^1$ by the ideal generated by
\[ I=\{\extd a_i\tens \extd b_i|\ a_i\extd b_i=0\}\subseteq\Omega^1
\tens_H\Omega^1.\]
Here $\extd$ extends to $\extd:\widetilde{\Omega^\cdot}\to
\widetilde{\Omega^\cdot}$ by
$\extd^2=0$ and the graded Leibniz rule
\[ \extd(\alpha\tens_H\beta)=\extd\alpha\tens_H \beta+(-1)^{|\alpha|}
\alpha\tens_H\extd\beta\]
for $\alpha$ of degree $|\alpha|$. Moreover,
$\widetilde{\Omega^\cdot}$ remains a super-Hopf algebra and $\extd$
commutes with its coproduct and antipode.
\end{propos}
\begin{proof}
We first observe that $I$ is itself a bicovariant
bimodule. For the coactions this follows from the fact that $\diff$ is
a bicomodule map. For the actions suppose $a_i\diff b_i=0$. Then
\[c \diff a_i\tens_H \diff b_i
 =\diff(c a_i)\tens_H \diff b_i-\diff c \tens_H a_i \diff b_i
 =\extd(ca_i)\tens_H\extd b_i\in I \]
for all $c\in H$, since $c a_i\diff b_i=0$. On the other side
\[\diff a_i\tens_H (\diff b_i) c=\diff a_i\tens_H \diff (b_i c)
 -\diff a_i\tens_H b_i \diff c
 =\diff a_i\tens_H \diff (b_i c)
 -\diff (a_i b_i)\tens_H \diff c + a_i \diff b_i \tens_H \diff c\in I\]
since $a_i\diff(b_i c)-a_i b_i\diff c=a_i(\diff b_i) c=0$.
We then define $\extd:\Omega^1\to
\Omega^1\tens_H\Omega^1/I=\Omega^2$ by
$\extd(a\extd b)=\extd a\tens_H\extd b$. This is well-defined
in virtue of the definition of $I$ and $\extd^2=0$, and also obeys the
Leibniz rule
$\extd (a\alpha)=\extd a\tens_H\alpha + a\extd\alpha$ (and similarly
with $a$ on the other side). We then extend to products of degree 1 forms by
the braided-Leibniz rule as stated. It is well-defined
on $\alpha\tens_H\beta$ in view of the above Leibniz rule with $a\in H$
and moreover itself obeys this Leibniz rule with $\alpha$ now of higher
degree. It maps to $\Omega^3=\Omega^2\tens_H\Omega^1+
\Omega^1\tens_H\Omega^2$.
In this way one constructs $\Omega^n$ and $\extd$ by induction.
In all higher
degrees $\Omega^n$ is the quotient of
$\Omega^1\tens_H\cdots\tens_H\Omega^1$ by
$\sum \Omega^1\tens_H\cdots\tens_H I\tens_H\cdots\tens_H\Omega^1$ as stated.
It is then straightforward to verify by induction that $\extd^2=0$ to all
orders. Next, we check that
$\widetilde{\Omega^\cdot}$ remains a super-Hopf algebra as a quotient
of $T_{-1}\Omega^1$.
There are general arguments for this, however in our case it is enough to
verify that $\cop I\subseteq I\tensor H + H\tensor I$ in
$T_{-1}\Omega^1$. Thus,
\begin{align*}
\cop(\diff a_i\tensor_H \diff b_i)
&=(\cop\diff a_i)\tensor_H (\cop\diff b_i)\\
&=(\diff\cop a_i)\tensor_H (\diff\cop b_i)\\
&=(\diff a_{i(1)}\tensor a_{i(2)} + a_{i(1)}\tensor \diff a_{i(2)})
  \tensor_H (\diff b_{i(1)}\tensor b_{i(2)}
   + b_{i(1)}\tensor \diff b_{i(2)})\\
&=\diff a_{i(1)}\tensor_H \diff b_{i(1)}
  \tensor a_{i(2)} b_{i(2)}
   + (\diff a_{i(1)}) b_{i(1)}
  \tensor a_{i(2)} \diff b_{i(2)}\\
 &\quad -a_{i(1)} \diff b_{i(1)}
  \tensor (\diff a_{i(2)}) b_{i(2)}
   + a_{i(1)} b_{i(1)}
  \tensor \diff a_{i(2)}\tensor_H \diff b_{i(2)}\\
&=\diff a_{i(1)}\tensor_H \diff b_{i(1)}\tensor a_{i(2)} b_{i(2)}
   + \diff (a_{i(1)} b_{i(1)})\tensor a_{i(2)} \diff b_{i(2)}\\
 &\quad -a_{i(1)} \diff b_{i(1)}\tensor \diff(a_{i(2)} b_{i(2)})
   + a_{i(1)} b_{i(1)}\tensor \diff a_{i(2)}\tensor_H \diff b_{i(2)}\\
&=(\diff a_i\tensor_H \diff b_i)\ro
  \tensor (\diff a_i\tensor_H \diff b_i)\rt
   + \diff (a_i \diff b_i)\lo\tensor (a_i \diff b_i)\lt\\
 &\quad -(a_i \diff b_i)\ro\tensor \diff (a_i \diff b_i)\rt
   +(\diff a_i\tensor_H \diff b_i)\lo
  \tensor (\diff a_i\tensor_H \diff b_i)\lt.
\end{align*}
While the second and third term are obviously zero, the first and the
fourth term fulfil the condition since $I$ is a bicomodule.
Similarly, one checks that $S I\subseteq I$.
Finally,
it is clear by induction that $\extd$ commutes with $\Delta$ and $S$
in view of the first order $\extd$ a bicomodule map. This step is the same
as in \cite{Brz:rem} for the Woronowicz-type exterior algebras.
\end{proof}

The maximal prolongation here is a quadratic algebra (over $H$) with relations
$I$ in degree 2. It is possible to generalise the construction to
general $\widetilde{\Omega^\cdot}_q$ similarly as a quotient of
$T_q\Omega^1$, but it
will no longer be quadratic. For example if $q$ is a primitive $n$-th root
of unity. Thus, for $n=3$ one must specify
\[ \extd:H\to\Omega^1,\quad \extd^2:H\to \Omega^1\tens_H\Omega^1,\]
where $\extd$ obeys the Leibniz rule (a given first order calculus) and
$\extd^2$ obeys the higher Leibniz rule
\[ \extd^2(ab)=(\extd^2 a)b+[2]_q\extd a\tens_H\extd b+a\extd^2 b\]
with $[2]_q=1+q$. We then define $\extd(a\extd b)=\extd
a\tens_H\extd b +a\extd^2b$ with a modified definition
\[ I=\{\extd
a_i\tens_H\extd b_i +a_i\extd^2 b_i|\ a_i\extd b_i=0\}\] in degree
2. Similarly, we define $\extd$ on degree 2 and higher by the
graded-Leibnitz rule with $q^{|\ |}$ in place of $(-1)^{|\ |}$ but must
now quotient further in degree 3 for this to be well-defined (a
cubic relation).  By construction, $\extd^2=\extd\circ
\extd$ and $\extd^3=0$. In this way one may build up a generalised complex
in the spirit of \cite{Ker:gra}.

One may similarly define a general exterior algebra $\Omega^\cdot$
associated to a first order bicovariant calculus $\Omega^1$ as any
super-Hopf algebra quotient of $T_{-1}\Omega^1$ in
Proposition~\ref{prop:tensor_algebra} with $H,\Omega^1$ in degrees
0,1 and such that $\extd$ extends as a differential graded algebra
and commutes with $\Delta,S$.

\begin{propos}
\label{prop:twist_diff}
Let $H$ be a Hopf algebra, $\chi:H\tensor H\to k$ a unital
2-cocycle. Then first order bicovariant differential calculi and
exterior (super-Hopf) algebras over $H$ and $H^\chi$ are in
one-to-one correspondence by the functor $\cgu$. $\cgu$ is trivial
on $\diff$.
\end{propos}
\begin{proof}
The degree 1 part as a bicovariant bimodule twists by
the functor
$\mathcal{G}^\chi$, as does the entire tensor super-Hopf algebra
$T_{-1}\Omega^1$ in Corollary~\ref{cor:com_twist_star}. By the same
arguments (Lemma~\ref{lem:ptwist_free}
and Theorem~\ref{thm:ptwist_bibimod})
this descends to any super-Hopf algebra quotient
$\Omega^\cdot$ of $T_{-1}\Omega^1$ and yields
$(\Omega^\cdot)^\chi$ as a
super-Hopf algebra of the desired form.
To see that this is compatible with the same
$\extd$ requires us to check the Leibnitz rule. This is,
\begin{align*}
\diff (\alpha\bullet \beta) &=\chi(\alpha\o \tensor \beta\o )\,
  \diff(\alpha\bt \beta\bt)
  \,\chi^{-1}(\alpha\th \tensor \beta\th )\\
&=\chi(\alpha\o \tensor \beta\o )\,
  ((\diff \alpha\bt) \beta\bt
  + (-1)^{|\alpha\bt|} \alpha\bt \diff \beta\bt)
  \,\chi^{-1}(\alpha\th \tensor \beta\th )\\
&=\chi((\diff \alpha)\o \tensor \beta\o )\,
  (\diff \alpha)\bt \beta\bt
  \,\chi^{-1}((\diff \alpha)\th \tensor \beta\th )\\
&\quad +(-1)^{|\alpha|} \chi(\alpha\o \tensor (\diff \beta)\o )\,
  \alpha\bt (\diff \beta)\bt
  \,\chi^{-1}(\alpha\th \tensor (\diff \beta)\th )\\
&=(\diff \alpha)\bullet \beta + (-1)^{|\alpha|} \alpha\bullet\diff \beta.
\end{align*}
Since the bicomodule and coproduct structure is unaffected by the
twist,  $\diff$
continues to be a bicomodule map and commute with the coaction
(so in particular,
$(\Omega^{1\chi},\extd)$ is a first order calculus in $H^\chi$).
Finally, the commutation of the antipode with $\diff$ is
\begin{align*}
\antip^\chi \diff \alpha &=U((\diff \alpha)\o ) \antip((\diff \alpha)\bt)
 U^{-1}((\diff \alpha)\th )\\
&=U(\alpha\o ) \antip(\diff \alpha\t ) U^{-1}(\alpha\th )\\ &=U(\alpha\o )
\diff\antip(\alpha\t ) U^{-1}(\alpha\th )\\ &=\diff\antip^\chi \alpha.
\end{align*}
\end{proof}

This tells us that given an exterior algebra $\Omega^\cdot$ on $H$, its twist
by $\chi$ is some other exterior algebra $\Omega^{\cdot\chi}$ on $H^\chi$.

\begin{corol}
The maximal prolongation is stable under twisting, i.e.
$\left(\widetilde{\Omega^\cdot}\right)^\chi$ is isomorphic  via $c^\chi$ to the
maximal prolongation of $\Omega^{1\chi}$.
\end{corol}
\begin{proof}
By Proposition~\ref{prop:twist_diff}, $(\widetilde{\Omega^\cdot})^\chi$ is an
exterior algebra. Since $\widetilde{\Omega^\cdot}$ is a quotient of
$T_{-1}\Omega^1$, we just have to ensure
that $\cnatu$, which identifies $(T_{-1}\Omega^1)^\chi$ with
$T_{-1}({\Omega^1}^\chi)$
according to Corollary~\ref{cor:com_twist_star}
identifies the corresponding ideal $I$ of
Proposition~\ref{prop:max_prolong} with its twisted counterpart. This is
\begin{align*}
I^\chi &=\{\diff a_i\tenshu \diff b_i | a_i\bullet \diff b_i=0\}\\
&=\{\diff a_i\tenshu \diff b_i |
 \chi(a_i\o\tensor b_i\o)\,
 a_i\t \diff b_i\t\,\chi^{-1}(a_i\th\tensor b_i\th)=0\}\\
&=\{\chi^{-1}(a_i\o\tensor b_i\o)\,
 \diff a_i\t\tenshu \diff b_i\t
 \,\chi(a_i\th\tensor b_i\th)
 | a_i \diff b_i=0\}\\
&=\cnatu\left(\{\diff a_i\tensh \diff b_i | a_i \diff b_i=0\}\right)
 =\cnatu(I).
\end{align*}
The third equality here is by substitution of
$\chi^{-1}(a_i\o\tensor b_i\o)\,
 a_i\t \tensor b_i\t\,\chi(a_i\th\tensor b_i\th)$
for $a_i\tensor b_i$.
\end{proof}

We also have a similar result for the Woronowicz exterior algebra
$\widehat{\Omega^\cdot}$. As explained in the Preliminaries this is also a
quotient of $T_{-1}\Omega^1$ but this time by skew-braid relations.

\begin{corol}
The Woronowicz exterior algebra is stable under
twisting, i.e.\ $\left(\widehat{\Omega^\cdot}\right)^\chi$ is isomorphic via $c^\chi$ to
the Woronowicz construction based on $\Omega^{1\chi}$.
\end{corol}
\begin{proof}
By Proposition~\ref{prop:twist_diff}, $(\widehat{\Omega^\cdot})^\chi$ is an
exterior algebra. Similarly to the preceding case, we view $\widehat{\Omega^\cdot}$
as a quotient of $T_{-1}\Omega^1$ and have to ensure
that $\cnatu$ is an intertwiner for the Woronowicz ideal
by which we quotient. But this ideal is given as the kernel of a linear
combination of identities and (bicovariant bimodule) braidings (see
Preliminaries) and $\cnat$ is an intertwiner for the braiding, so this
is satisfied.
\end{proof}

Finally, we have explained in the Preliminaries that bicovariant
calculi are of the form $\Omega^1=V\tens H$ where $V$ is a quotient of
$\ker\eps\in\catlc{H}$. Here $\ker\eps$ and $H$ itself are crossed modules
by left multiplication and the left adjoint coaction. To complete our picture,

\begin{propos}
\label{prop:alpha}
We denote by $\pi:\ker\eps\to V$ the projection corresponding
to $\Omega^1$ on $H$. Then the calculus $\Omega^{1\chi}$ on
$H^\chi$ corresponds to the projection
$\pi^\chi=\pi\circ\alpha^{-1}$, where
\[ \alpha:\cfunct^\chi(H)\to H^\chi,\qquad
 \alpha(h)=\chi^{-1}(h\o \tensor\antip h\fiv )\,
  \chi(h\th \tensor\antip h\fo )\, h\t \]
is an isomorphism in $\catlc{H^\chi}$ with inverse
\[ \alpha^{-1}(h)=\chi^{-1}(h\o \antip h\th \tensor
 h\fo )\, h\t \]
and restricts to an isomorphism
${\cfunct}^\chi(\ker\eps)\isom\ker\eps$.
\end{propos}
\begin{proof}
We deduce this from our theory of twisting of first order
differential calculi. Using Proposition~\ref{prop:twist_diff}, we
know that $(V\tensor H)^\chi=\cgu(V\tensor H)$ is a differential
calculus over $H^\chi$ and hence of the form $V^\chi\tensor H^\chi$
for some $V^\chi$ a quotient of $\ker\eps\subset H^\chi$. This
$V^\chi$ is the right-invariant subspace of the corresponding
differential calculus. On the other hand, the projection for the
subspace $V$ can be obtained explicitly from its associated
calculus by $\pi(a)=\diff a\o
\cdot\antip a \t $, and similarly for $\pi^\chi$ using the Hopf
algebra $H^\chi$. Thus, from Proposition~\ref{prop:twist_diff}, we
find,
\begin{align*}
\pi^\chi(a) &=\diff a\o \bullet\antip^\chi a\t \\
 &=\diff a\o \bullet\antip a\th \,U(a\t ) U^{-1}(a\fo )\\
 &=\chi(a\o \tensor \antip a\sev )\,
   \diff a\t \cdot\antip a\six \,
   \chi^{-1}(a\th \tensor \antip a\fiv )\,
   U(a\fo ) U^{-1}(a\eig)\\
 &=\chi(a\o \tensor \antip a\fo )\,
   \diff a\t \cdot\antip a\th \,
   U^{-1}(a\fiv )\\
 &=\chi^{-1}(a\o \antip a\fo \tensor a\fiv )\,
   \diff a\t \cdot\antip a\th \\
 &=\pi(\chi^{-1}(a\o \antip a\th \tensor a\fo )\,
   a\t )
\end{align*}
as stated. In particular, we apply these arguments to the universal
differential calculus, which corresponds to $V=\ker\cou\subset H$
and hence obtain $\alpha^{-1}$ as a linear isomorphism on
$\ker\eps$. By construction it must in fact identify
$\cfunct^\chi(\ker\eps)$ as a crossed module in
$\catlc{H^\chi}$ with $\ker\eps\subset H^\chi$. Finally, it is
trivial to check that it extends to the whole of $H=\ker\cou\oplus
k 1$ with $\alpha(1)=1$ by the same formula.
\end{proof}

This tells us also that if $V=\ker\eps/M$
by some $\Ad$-stable ideal $M$ then the corresponding ideal for
$V^\chi$ is $M^\chi=\alpha\circ\cfunct^\chi(M)$. Note that the
isomorphism of
crossed modules $\alpha:\cfunct(H)\to H^\chi$ is somewhat nontrivial to obtain
by normal Hopf algebraic methods; the intertwining of the adjoint coactions
alone is a first step introduced (in a dual form) in \cite{GurMa:bra}, which
was in fact the starting point behind the present paper.

We turn now to some general applications of our twisting theory, beyond the
Planck scale Hopf algebra to be studied in the next section. In fact, many
interesting quantum groups are related by twisting and our theory allows the
construction of their quantum differential calculi one from the other. We give
two important general constructions where the quantum group is a twisting of a
tensor product quantum group. But for tensor product Hopf algebras it is easy
to obtain calculi given calculi on the pieces. Indeed, if
$(H_1,\Omega_1^\cdot)$
and $(H_2,\Omega_2^\cdot)$ are two Hopf algebras equipped with
calculi and associated exterior super-Hopf algebras, it is clear that
 $\Omega^\cdot\equiv \Omega_1^\cdot\und\tensor\Omega_2^\cdot$ (the super
 tensor product) is an exterior super-Hopf algebra  over
$H:=H_1\tensor H_2$.
$\diff:\Omega^\cdot\to\Omega^\cdot$ is given by the Leibnitz rule.
Restriction to degree 1 yields the corresponding construction for
first order bicovariant calculi $\Omega^1$, i.e. $\Omega^1=\Omega^1_1
\tensor H_2\oplus H_1\tensor\Omega^1_2$ given $\Omega^1_i$, with $\extd$
defined via the Leibnitz rule.

\begin{propos}
Given first order calculi $\diff:H\to\Omega^1$ on a
finite-dimensional Hopf algebra $H$ and $\diff:H^{*\rm
op}\to\Omega^1_{*\rm op}$ on $H^{*\rm op}$, we obtain a calculus on
the quantum double $D(H)$ by
\[ \Omega_{D(H)}^1=\Omega^1\tensor H^{*\rm op}\oplus
H\tensor\Omega^1_{*\rm op}\]
and the additional bimodule structure
\[h\, \diff\phi =\<\antip h\o ,\phi\o \>
 (\diff\phi\t)  h\t  \<h\th ,\phi\th \>,\quad
 \phi\, \diff h =\<h\o ,\phi\o \>
 (\diff h\t)  \phi\t  \<\antip h\th ,\phi\th \>\]
for all $h\in H$ and $\phi\in H^*$. Here $\extd$ restricts to the
given one on $H,H^{*\rm op}\subseteq D(H)$.
\end{propos}
\begin{proof}
It is known that the quantum double $D(H)=H\dcross H^{*\rm op}$ is a cotwist
of $H\tens H^{*\rm op}$ by the cocycle
\[ \chi((h\tens\phi)\tens (g\tens\psi))=\eps(h)\<\phi,Sg\>\eps(\psi).\]
See \cite{Ma:book}. We then apply Proposition~\ref{prop:twist_diff},
i.e.\ the functor
$\mathcal{G}^\chi$. The vector space, coactions and $\extd$ are not changed
from the tensor product calculus under the functor, but the bimodule
structures are,
as shown. These are easily computed from the form of the cocycle.
\end{proof}

This is a different and rather more geometrical approach to the construction
of the differential calculi on $D(H)$ than the one in \cite{Ma:cla} based on
its representation theory. In a similar spirit we may consider the general
double cross product Hopf algebra  $H\dcross_{\CR} H$ associated to any
dual-quasitriangular Hopf algebra $H$ \cite{Ma:book}.

\begin{propos} Let $H$ be dual-quasitriangular. Given two calculi
$\extd:H\to \Omega^1_L$ and $\extd:H\to \Omega^1_R$ we obtain a calculus
on $H\dcross_{\CR} H$ by
\[ \Omega^1=\Omega^1_L\tens H\oplus H\tens\Omega^1_R\]
and the additional bimodule structure
\[
h\, \diff a =\CR^{-1}(h\o, a\o)
 (\diff a\t)  h\t  \CR(h\th ,a\th ),\quad
a\, \diff h =\CR(h\o , a\o )
 (\diff h\t)  a\t  \CR^{-1}(h\th ,a\th )
\]
for all $h$ in the first copy of $H$ and $a$ in the second. Here
$\extd$ restricted to each copy of $H$ is the given one.
\end{propos}
\begin{proof}
We proceed in the same fashion as above. The cocycle is given by
\[ \chi((h\tens a)\tens(g\tens b))=\eps(h)\CR^{-1}(a\tens g)\eps(b),
\quad\forall h,g,a,b\in H\]
and allows us to write $H\dcross_{\CR}H=(H\tens H)^\chi$. See \cite{Ma:book}.
Here the two copies of $H$ remain sub-Hopf algebras as for the quantum
double case
above, and the remaining computation is similar. \end{proof}

In particular, a standard formulation of the q-Lorentz group is as
$SU_q(2)\dcross_{\CR}SU_q(2)$, a twist of the q-Euclidean rotation group as
$SU_q(2)\tens SU_q(2)$. Usually the calculi on these are obtained and
studied separately, but the above proposition constructs one from the other.
This extends the `quantum Wick rotation' in
\cite{Ma:euc} to the construction of bicovariant
differential calculi. In both the above constructions the factors
appear as sub(Hopf) algebras and in this case there is a well-defined
notion of a calculus being decomposable or built up
from calculi on the factors; the above results fully classify such decomposable
calculi.

\section{Planck scale Hopf algebra as a cotwist and its differential
geometry}

The remainder of the paper applies the preceding results to one
particular Hopf algebra, namely the bicrossproduct Hopf algebra
 $\C[x]\bicross_{\hbar,\grav}\C[p]$ introduced in \cite{Ma:pla}. For
the purposes of the present section we work algebraically with
$g=e^{-\frac{x}{\grav}}$ and $g^{-1}$ instead of $x$. Then the
explicit formulae (as stated in the introduction) are
\begin{gather*}
[p,g]=\imath A(1-g)g,\\
\cop p = p\tensor g + 1\tensor p,\qquad
\cop g = g\tensor g,\\
\antip p = -p g^{-1},\quad  \cou p = 0,\qquad
\antip g = g^{-1},\quad  \cou g = 1,
\end{gather*}
where $A=\frac{\hbar}{\grav}$. Also, as $\hbar\to 0$ (corresponding
to $A\to 0$), we obtain $C(B_+)$ (in an algebraic form) where
\[
 B_+=\R\rcross_\grav \R,\quad
 (x,p)(x',p')=(x+x',pe^{-\frac{x'}{\grav}}+p')
\]
as explained in \cite{Ma:pla}. We consider $B_+$ to be the classical
phase space underlying the quantum system described by the
Plank-scale Hopf algebra. In terms of the exponentiated coordinate
in place of $x$, the group law of $B_+$ is that of a matrix group
\[
 \begin{pmatrix} g & 0\\ p & 1 \end{pmatrix}
 \begin{pmatrix} g'& 0\\ p'& 1 \end{pmatrix}
=\begin{pmatrix} gg'& 0\\ pg'+p' & 1 \end{pmatrix}.
\]
The notation $\pla$ reflects the construction of this Hopf algebra
as a bicrossproduct, i.e.\ a semidirect product as an algebra and a
semidirect coproduct as a coalgebra by certain actions and
coactions arising from a Lie group
factorisation \cite{Ma:hop}\cite{Ma:mat}\cite{Ma:pla}. This aspect will be used
extensively in Section~5.

\subsection{The cocycle and differential calculi}

Our starting points are the known facts that the Hopf algebra
$\pla$ is of a self-dual form and at the same time a
twisting (of the coproduct) of $U(\cb_+)$ where $\cb_+$ is the
Borel subalgebra of $sl_2$. Combining these observations, one may
expect that it is also a product twist by a cocycle. This turns out
to be the case.

\begin{propos}
\label{prop:pla_cocycle}
$\chi$ defined by
\begin{gather*}
\chi = (\cou\tensor\cou)\circ
 \exp(\imath A \pdiff{p}\tensor\pdiff{g^{-1}}),\qquad
\chi^{-1} = (\cou\tensor\cou)\circ
 \exp(\imath A \pdiff{p}\tensor\pdiff{g})
\end{gather*}
is a unital 2-cocycle on $C(B_+)$, and $\pla=C(B_+)^\chi$.
\end{propos}
\begin{proof}
In order to show that $\chi$ is a unital 2-cocycle we have to show its
invertibility, the cocycle condition and the unitality (see
Preliminaries). It will be
useful to have the explicit expressions of $\chi$ and $\chi^{-1}$ on
a basis $\{p^n g^r | n\in\N_0, r\in\Z\}$ of $\cbp$:
\begin{gather*}
\chi(p^n g^r\tensor p^m g^s)
 =\delta_{m,0} (\imath A)^n \prod_{k=0}^{n-1} (-s-k),\qquad
\chi^{-1}(p^n g^r\tensor p^m g^s)
 =\delta_{m,0} (\imath A)^n \prod_{k=0}^{n-1} (s-k).
\end{gather*}
For the invertibility we require
\begin{gather*}
\chi(a\o\tensor b\o)\,\chi^{-1}(a\t\tensor b\t)
 =\cou(a)\cou(b),\qquad
\chi^{-1}(a\o\tensor b\o)\,\chi(a\t\tensor b\t)
 =\cou(a)\cou(b).
\end{gather*}
To see this we take
$a=p^n g^r$ and $b=p^m g^s$ and the first expression becomes
\begin{equation*}
\begin{split}
&\sum_{k,l} \binom{n}{k} \binom{m}{l}
 \chi(p^k g^r\tensor p^l g^s)\,
 \chi^{-1}(p^{n-k} g^{k+r}\tensor p^{m-l} g^{l+s})\\
&=\sum_{k,l} \binom{n}{k} \binom{m}{l}
 \delta_{l,0}\, \delta_{m,l} (\imath A)^n
 \prod_{i=0}^{k-1} (-s-i) \prod_{j=0}^{n-k} (l+s-j)\\
&=\delta_{m,0} (\imath A)^n \sum_{k} \binom{n}{k}
 \prod_{i=0}^{k-1} (-s-i) \prod_{j=0}^{n-k} (s-j)\\
&=\delta_{m,0}\, \delta_{n,0}=\cou(p^n g^r) \cou(p^m g^s).
\end{split}
\end{equation*}
We have used
\[\sum_{k=0}^n \binom{n}{k}
 \prod_{i=0}^{k-1} (-s-i) \prod_{j=0}^{n-k} (s-j)=\delta_{n,0},
 \quad\forall s\in\Z\ \forall n\in\N_0.\]
While this is obvious for $s=0$ it follows easily by induction
for $s\neq 0$.
Note that the exchange of $\chi$ and $\chi^{-1}$ in the above
calculation is equivalent to
replacing $s$ by $-s$. Thus follows the second equation as well.
Next, for the cocycle condition
we take $p^n g^r$, $p^m g^s$, $p^l g^t$ for $h$, $g$, $f$ in
(\ref{cochi}).
The left hand side evaluates to
\begin{equation*}
\begin{split}
&\sum_{j,k} \binom{m}{j} \binom{l}{k}
 \chi(p^j g^s\tensor p^k g^t)\,
 \chi(p^n g^r\tensor p^{m-j+l-k}g^{j+s+k+t})\\
&=\sum_{j,k} \binom{m}{j} \binom{l}{k}
 \delta_{k,0} \delta_{m-j+l-k,0} (\imath A)^{j+n}
 \prod_{i=0}^{j-1}(-t-i) \prod_{h=0}^{n-1} (-j-s-k-t-h)\\
&=\sum_{j} \binom{m}{j}
 \delta_{m-j+l,0} (\imath A)^{j+n}
 \prod_{i=0}^{j-1}(-t-i) \prod_{h=0}^{n-1} (-j-s-t-h)\\
&=\delta_{l,0} (\imath A)^{n+m}
 \prod_{i=0}^{m-1}(-t-i) \prod_{h=0}^{n-1} (-m-s-t-h).
\end{split}
\end{equation*}
The right hand side is
\begin{equation*}
\begin{split}
&\sum_{j,k} \binom{n}{j} \binom{m}{k}
 \chi(p^j g^r\tensor p^k g^s)\,
 \chi(p^{n-j+m-k}g^{j+r+k+s}\tensor p^l g^t)\\
&=\sum_{j,k} \binom{n}{j} \binom{m}{k}
 \delta_{k,0} \delta_{l,0} (\imath A)^{n+m-k}
 \prod_{i=0}^{j-1}(-s-i) \prod_{h=0}^{n-j+m-k-1} (-t-h)\\
&=\delta_{l,0} (\imath A)^{n+m} \sum_{j} \binom{n}{j}
 \prod_{i=0}^{j-1}(-s-i) \prod_{h=0}^{n-j+m-1} (-t-h)\\
&=\delta_{l,0} (\imath A)^{n+m}
 \prod_{i=0}^{m-1}(-t-i) \prod_{h=0}^{n-1} (-m-s-t-h).
\end{split}
\end{equation*}
The last equality can be easily checked by induction in $n$.
Finally, the unitality in (\ref{cochi}) follows easily from the
explicit formula for $\chi$.
It remains to check that the twist of $\cbp$ defined by $\chi$ is
indeed $\pla$. For that, it is sufficient to check the commutator
between $p$ and $g$.
For clarity, we distinguish the twisted product from the untwisted one
by denoting the former with a $\bullet$.
\begin{equation*}
\begin{split}
g\bullet g &= \chi(g\tensor g)\, g g\, \chi^{-1}(g\tensor g)
 = g g\\
p\bullet g &= \chi(p\tensor g)\, g g\, \chi^{-1}(g\tensor g)
 +\chi(1\tensor g)\, p g\, \chi^{-1}(g\tensor g)
 +\chi(1\tensor g)\, g\, \chi^{-1}(p\tensor g)\\
 &= - \imath A g g + p g + \imath A g\\
g\bullet p &= \chi(g\tensor p)\, g g\, \chi^{-1}(g\tensor g)
 +\chi(g\tensor 1)\, g p\, \chi^{-1}(g\tensor g)
 +\chi(g\tensor 1)\, g\, \chi^{-1}(g\tensor p) = g p
\end{split}
\end{equation*}
In particular we obtain
\[p\bullet g - g \bullet p = \imath A (1-g) g = \imath A (1-g) \bullet g ,\]
which is the correct relation in $\pla$.
\end{proof}

Let us now turn to the differential calculi on $\pla$. We first
repeat a result from \cite{Oe:cla} about the untwisted Hopf algebra
$\cbp$.

\begin{propos}
\cite{Oe:cla}
\label{prop:cbp_cla}
(a) Finite dimensional differential calculi $\Omega^1$ on \cbp{}
are in one-to-one correspondence to non-empty finite sets
$I\subset\N$ and have dimension $(\sum_{n\in I} n)-1$.

(b) The differential calculus of dimension $n\ge 2$
corresponding to $\{n,1\}$ has a right invariant basis
$\eta_0,\dots\eta_{n-1}$ so that
\begin{gather*}
\diff g=g \eta_0,\qquad \diff g=g\eta_1,\\
[g,\eta_k]=0,\qquad
[p,\eta_k]=
 \begin{cases}
  0& \text{if $k=0$ or $k=n-1$}\\
  g \eta_{k+1}& \text{if $0<k<n-1$}
 \end{cases}\\
\beta_L(\eta_k)=
 \begin{cases}
  g^{-k}\tensor \eta_k& \text{if}\ k\neq 1\\
  g^{-1}\tensor \eta_1+g^{-1}p\tensor \eta_0& \text{if}\ k=1
 \end{cases}
\end{gather*}

(c) The differential calculus of dimension $n-1\ge 1$
corresponding to $\{n\}$ is the same as (b) except that $\eta_0=0$.
\end{propos}
\begin{proof}
We refer to \cite[Prop.~3.5]{Oe:cla}. Note however, that in
\cite{Oe:cla} different generators were used for $\cbp$, so that
the coproduct appeared in a different form. In our conventions the
crossed submodule $M\subset\ker\cou$ corresponding to $\{n,1\}$ is
generated by $(g-1)(g-1),p(g-1),\dots,p^n$ as a crossed module.
Denoting the equivalence classes of $g-1,p,\dots,p^{n-1}$ in
$\ker\cou/M$ by $\eta_0,\dots,\eta_{n-1}$, we obtain the derivative
and commutation relations as stated. For the left coaction note
that the left adjoint action on $\cbp$ takes the form
\[\adl(f(g)p^k)=\sum_{t=0}^k \binom{k}{t} g^{-k}p^{k-t}\tensor
 f(g)(g-1)^{k-t}p^t.\]
$\beta$ is then obtained by composition with the projection to
$\ker\cou/M$.
\end{proof}

We can now apply our twisting theory in Sections~2,3 to solve the
classification problem for calculi on the Planck scale Hopf
algebra,

\begin{propos}
\label{prop:pla_cla}
(a) Finite dimensional differential calculi $\Omega^1$ on $\pla$
are in one-to-one correspondence to non-empty finite subsets
$I\subset\N$ with dimensions as in Proposition~\ref{prop:cbp_cla}.

(b) The differential calculus of dimension $n\ge 2$
corresponding to $\{n,1\}$ has a right invariant basis
$\eta_0,\dots\eta_{n-1}$ so that
\begin{gather*}
\diff g=g \eta_0,\qquad \diff p=g\eta_1,\\
[g,\eta_k]=
 \begin{cases}
  0 & \text{if}\ k\neq 1\\
  \imath A g \eta_0 & \text{if}\ k=1
 \end{cases}
\qquad [p,\eta_k]=\imath A k g\eta_k +
 \begin{cases}
  0 & \text{if $k=0$ or $k=n-1$}\\
  g \eta_{k+1}& \text{if $0<k<n-1$}
 \end{cases}
\end{gather*}
and $\beta$ of the same form as in Proposition~\ref{prop:cbp_cla}.

 (c) The
differential calculus of dimension $n-1\ge 1$ corresponding to
$\{n\}$ is the same as (b) except that $\eta_0=0$.
\end{propos}
\begin{proof}
We apply Proposition~\ref{prop:twist_diff} to Proposition~\ref{prop:cbp_cla}.
Part (a) remains unchanged. For part (b) we calculate the twisted
actions in terms of the untwisted ones (using a $\bullet$ to denote the
twisted ones).
\begin{align*}
g\bullet\eta_k &=\chi(g\tensor g^{-k})\, g\eta_k+\delta_{k,1}
 \chi(g\tensor g^{-1}p)\, g\eta_0 = g\eta_k\\
\eta_k\bullet g &=\chi(g^{-k}\tensor g)\, \eta_k g
 +\delta_{k,1}\chi(g^{-1}p\tensor g)\,\eta_0 g
 = \eta_k g - \imath A\delta_{k,1}\eta_0 g\\
p\bullet\eta_k &=\chi(p\tensor g^{-k})\, g\eta_k
 +\chi(1\tensor g^{-k})\, p\eta_k
 +\delta_{k,1}(\chi(p\tensor g^{-1}p)\, g\eta_0
 +\chi(1\tensor g^{-1}p)\, p\eta_0)\\
 &= \imath A k g\eta_k + p \eta_k\\
\eta_k\bullet p &=\chi(g^{-k}\tensor p)\, \eta_k g
 +\chi(g^{-k}\tensor 1)\, \eta_k p
 +\delta_{k,1}(\chi(g^{-1}p\tensor p)\, \eta_0 g
 +\chi(g^{-1}p\tensor 1)\, \eta_0 p)
 = \eta_k p
\end{align*}
This gives the new commutators and
the expressions for the differentials. For the
coaction we observe that $g^{-1}\bullet p=g^{-1} p$ so that its form
does not change. Part (c) remains unchanged.
\end{proof}

For the remainder of the section we concentrate on the
calculus $\{2,1\}$ which is the quantisation of the standard
classical calculus on $B_+$. We can use the twisting theory to
quantise in fact the entire exterior algebra in this case.

\begin{propos}
The exterior algebra $\Omega^\cdot$ of $\pla$ corresponding via
twisting to the
classical one of $\cbp$ has the following properties.
The first order calculus
has a basis $\{\xi,\eta\}$ of right-invariant
1-forms with
\begin{gather*}
\diff g = g \xi, \qquad \diff p = g \eta, \qquad
[a,\xi]=0,\qquad [a,\eta]=\imath A \diff a, \qquad
 \forall a\in\pla,\\
\beta_L(\xi)=1\tens\xi,\qquad \beta_L(\eta)=g^{-1}\tensor
\eta+g^{-1}p\tensor\xi.
\end{gather*}
The 2-forms have relations
\begin{gather*}
\xi\wedge\xi=0,\qquad \eta\wedge\xi=-\xi\wedge\eta,
 \qquad \eta\wedge\eta=\imath A \xi\wedge\eta,\\
\diff\xi=0,\qquad
 \diff\eta=\eta\wedge\xi.
\end{gather*}
As a $\Z_2$-graded Hopf algebra $\Omega^\cdot$ has the structure
\begin{gather*}
\cop\xi=\xi\tensor 1+1\tensor\xi,\quad
\cop\eta=g^{-1}\tensor\eta+g^{-1}p\tensor\xi-\eta\tensor 1,\\
\cou(\xi)=\cou(\eta)=0,\qquad
\antip\xi=-\xi,\quad\antip\eta=-g\eta+p\xi.
\end{gather*}
\end{propos}
\begin{proof}
For the first order calculus we define $\xi\equiv \eta_0$ and
$\eta\equiv \eta_1$. The commutation relations in the $n=2$ case
(b) of Proposition~\ref{prop:pla_cla} become as stated. Next, the
classical space
of 2-forms on $\cbp$ is spanned by $\xi\wedge\eta=-\eta\wedge\xi$.
Denoting the wedge product on $\pla$ by $\wedge_\bullet$ we have
\begin{align*}
\eta\wedge_\bullet\eta &=\chi(g^{-1}\tensor g^{-1})\,\eta\wedge\eta
 +\chi(g^{-1} p\tensor g^{-1})\,\xi\wedge\eta\\
 &\quad +\chi(g^{-1}\tensor g^{-1} p)\,\eta\wedge\xi
 +\chi(g^{-1} p\tensor g^{-1} p)\,\xi\wedge\xi\\
 &=\imath A \xi\wedge\eta.
\end{align*}
The other wedge products involving $\xi$ and $\eta$ are identical
to the classical ones due to the bi-invariance of $\xi$. This leads
to the relations stated. Finally, for the differentials of the
1-forms observe (the twisted and untwisted wedge products are the
same here)
\begin{align*}
\diff\xi &=\diff(g^{-1}\diff g)=\diff g^{-1}\wedge\diff g
 =-g^{-1}\xi\wedge g\xi=0,\\
\diff\eta &=\diff(g^{-1}\diff p)=\diff g^{-1}\wedge\diff p
 =-g^{-1}\xi\wedge g\eta=\eta\wedge\xi.
\end{align*}
The coproduct and antipode are readily obtained using
Proposition~\ref{prop:tensor_algebra}. The exterior algebra here
coincides with the Woronowicz prolongation of the first order part.
\end{proof}

In terms of generators $x$ and $p$ the exterior algebra is generated
by $\extd x,\extd p$ with the relations
\begin{gather*}
a\extd x=(\extd x)a,\quad a\extd p=(\extd p)a
 +\frac{\imath\hbar}{\grav}\extd a,\\
\extd x\wedge\extd x=0,\quad
 \extd x\wedge\extd p=-\extd p\wedge \extd x,\quad
\extd p\wedge\extd p=0.
\end{gather*}
From this we see explicitly that in the classical limit $\hbar\to 0$ we
obtain the usual exterior algebra on $B_+$. By contrast, the other
limit $\grav\to0$ is highly singular with these generators, so that
the exterior algebra is not even defined in this case. In other
words, the presence of `gravity' in the form of $\grav$ restores
the geometrical picture not visible in flat space quantum
mechanics.

We also know from the Preliminaries that associated to a first
order calculus is a quantum tangent space. The right-invariant
derivatives are generated by elements of $L$ dual to $V$ and obey a
braided Leibniz rule.

\begin{propos}
\label{prop:right_deriv}
Let $\{\xi^*,\eta^*\}$ be the basis of $L$ dual to the basis
$\{\xi,\eta\}$ above. Then,
\begin{align*}
\partial_\xi(:f(g,p):)&=
 :g\pdiff{g} f(g,p+\imath A)+g(f(g,p+\imath A)-f(g,p)): ,\\
\partial_\eta(:f(g,p):)&=:g\frac{f(g,p+\imath A)-f(g,p)}{\imath A}: .
\end{align*}
\end{propos}
\begin{proof}
We observe that $\diff g^n=(\diff g) n g^{n-1}$ and $\diff p^n=(\diff
p)\frac{(p+\imath A)^n-p^n}{\imath A}$ (this can be easily checked by
induction), so that
\begin{align*}
\diff (g^n p^m) &=(\diff g^n)p^m + g^n\diff p^m\\
 &=(\diff g^n)p^m + g^n (\diff p) \frac{(p+\imath A)^m-p^m}{\imath A}\\
 &=(\diff g^n)p^m + (\diff p)g^n \frac{(p+\imath A)^m-p^m}{\imath A}
   +\imath A (\diff g^n) \frac{(p+\imath A)^m-p^m}{\imath A}\\
 &=(\diff g^n)(p+\imath A)^m + (\diff p) g^n \frac{(p+\imath A)^m-p^m}
  {\imath A}\\
 &=\xi\, n g^n (p+\imath A)^m + (\eta\, g+ \imath A\, \xi\, g)
     g^n\frac{(p+\imath A)^m-p^m}{\imath A}\\
 &=\xi\, \left(n g^n (p+\imath A)^m + g^{n+1}
    \left((p+\imath A)^m-p^m\right)\right)
 +\eta\, g^{n+1} \frac{(p+\imath A)^m-p^m}{\imath A},
\end{align*}
which we compare with the property $\diff f=
\xi\,\partial_\xi(f)+\eta\,\partial_\eta(f)$ of the partial derivatives.
\end{proof}

In terms of coordinates $x$, $p$ we can similarly write the action of the
basis of $L$ dual to $\{\diff x,\eta\}$ as
\begin{align}
\label{rightderiv1}
\partial_x(:f(x,p):)&=
 :\pdiff{x} f(x,p+\frac{\imath \hbar}{\grav})
  - \frac{e^{-\frac{x}{\grav}}}{\grav}
  (f(x,p+\frac{\imath\hbar}{\grav})-f(x,p)): ,\\
\label{rightderiv2}
\partial_\eta(:f(x,p):)&=\frac{\grav}{\imath\hbar}:e^{-\frac{x}{\grav}}
 (f(x,p+\frac{\imath\hbar}{\grav})-f(x,p)):
\end{align}
(here $\del_x$ denotes the action of the basis element dual to
$\extd x$ by a slight abuse of notation).

Finally, for completeness we note that all these formulae are for
right-invariant differential forms. There is an equally good theory
based on $L,V\in\catrc{H}$ and left-invariant partial derivatives.
We take a left-invariant basis of the 1-forms to be
$\{\xi=g^{-1}\extd g,\bar\eta =\extd p-p g^{-1}\extd g\}$. The
relations of the the calculus become
\begin{gather*}
[a,\xi]=0, \qquad [a,\bar\eta ]=\imath A \diff a, \qquad
 \forall a\in\pla,\\
\beta_R(\xi)=\xi\tensor 1,\qquad
\beta_R(\bar\eta )=\bar\eta \tensor g-\xi\tensor p,\\
\xi\wedge\xi=0,\qquad \bar\eta \wedge\xi
 =-\xi\wedge\bar\eta, \qquad \bar\eta \wedge\bar\eta
 =\imath A \bar\eta \wedge\xi.
\end{gather*}
Moreover, the differential in $\Omega^\cdot$ is generated by
(graded) commutation with the element
$\theta=-\frac{1}{2}(\eta+\bar\eta )$ as
\eqn{theta}{ [\theta,\alpha]= \imath A\diff\alpha,\qquad
\forall\alpha\in\Omega^\cdot .}
This is a step towards a Connes spectral triple description of
this calculus, to be considered elsewhere.
The generator $\frac{1}{\imath A}\theta$ is singular in the limit $A\to 0$
($\hbar\to 0$)
so that the presence of $\hbar$ allows for nicer properties of the
differential calculus than visible classically. This is a typical feature
of q-deformation known for calculi on standard quantum groups.

We likewise bar the braided derivations in the left-invariant
version of the theory to avoid confusion with the above
right-invariant ones. The left-invariant derivations corresponding
to $\{\diff x,\bar\eta \}$ in the sense $\diff
f=(\bar\partial_{x}f)\diff x+(\bar\partial_{\eta}f)\bar\eta $ are
\eqn{leftdel1}{
\bar\partial_{x}(:f(x,p):)=:\pdiff{x}f(x,p)
 +\frac{p}{\imath\hbar}\,
 (f(x,p-\frac{\imath\hbar}{\grav})-f(x,p)):,}
\eqn{leftdel2}{\bar\partial_{\eta}(:f(x,p):)=-\frac{\grav}{\hbar}
:(f(x,p-\frac{\hbar}{\grav})
 -f(x,p)):.}

\subsection{Quantum Poisson bracket}

We conclude this section with some elements of `quantum Poisson
geometry'. We recall first of all the classical situation. In fact,
for any twisting of a commutative Hopf algebra $H$ by a cocycle
$\chi$ admitting a reasonable expansion in a parameter $\hbar$ (so
that the deformation is flat) one knows on general grounds that the
commutative Hopf algebra is (an algebraic version of) a Poisson-Lie
group. As for any flat deformation, the Poisson bracket is provided
by the $\hbar\to 0$ part of $\frac{\imath}{\hbar}[\ , ]$ (the
leading part of the commutator). When $H=\C[G]$ is an algebraic
group of Lie type one can say rather more. We let $\cg$ be the Lie
algebra of $G$. If the cocycle $\chi_\hbar$ has the form
\[ \chi_\hbar(a\tens b)=a(e)b(e)
 +\frac{\hbar}{\imath} \<\tilde{\sigma},\extd a\tens\extd b\>(e)
 +\O(\hbar^2)\]
where $\sigma\in\cg\tens\cg$ and $\tilde\sigma$ denotes
the extension as a left-invariant bivector field. Then it is known
from Drinfeld's theory \cite{Dri:hop} that
\[ r=\sigma-\tau(\sigma)\]
($\tau$ the flip map) is a triangular solution of the Classical Yang-Baxter
equation. Moreover,
\[ \{a,b\}=\<\tilde r-\tilde r^R,\extd a\tens \extd b\>\]
is the Poisson bracket of which $H^\chi$ is the quantisation, and
which makes $G$ a Poisson-Lie group. \cite{Dri:hop} proves in fact
the converse to this (the formal existence of $\chi_\hbar$) but the
above is also covered. In our case of the Planck scale Hopf
algebra, $G=B_+$ with Lie algebra $\cb_+$ spanned by $x_0,x_1$ with
relations $[x_0,x_1]=\frac{\imath}{\grav}x_1$ becomes a triangular Lie
bialgebra with
\[ r=x_0\tens x_1-x_1\tens x_0. \]
The resulting Lie bialgebra has the Lie bicross sum form $\bicross$
of which the bicrossproduct Hopf algebras are quantisations, see
\cite{Ma:ista}. Note that the corresponding Poisson bracket, as
with all (quasi)triangular Poisson-Lie groups, cannot be symplectic
since it must vanish at least at the group identity.

\begin{propos}
\label{prop:cbp_poisson}
The Poisson bracket on $\cbp$, for which the cocycle $\chi$ of
Proposition~\ref{prop:pla_cocycle} provides the quantisation, is
\[ \{a,b\}=(e^{-\frac{x}{\grav}}-1)\left( \pdiffr{a}{x}\pdiffr{b}{p}
-\pdiffr{b}{x}\pdiffr{a}{p}
\right).\]
\end{propos}
\begin{proof}
Expanding $\chi$ of Proposition~\ref{prop:pla_cocycle} in $\hbar$
and expressing everything in terms of the coordinates $x,p$ yields
\begin{align*}
a\bullet b &=a b +\imath\hbar \cou\left(\pdiff{p}a\o\right)
 \cou\left(e^{-\frac{x}{\grav}}\pdiff{x}b\o\right)a\t b\t\\
&\quad -\imath\hbar a\o b\o\cou\left(\pdiff{p}a\t\right)
 \cou\left(e^{\frac{x}{\grav}}\pdiff{x} b\t\right)
 +\O\left(\hbar^2\right)\\
&=a b +\imath\hbar \left(e^{-\frac{x}{\grav}}\pdiffr{a}{p}\pdiffr{b}{x}
 -\pdiffr{a}{p} \pdiffr{b}{x}\right)
 +\O\left(\hbar^2\right)\\
&=a b +\imath\hbar \left(e^{-\frac{x}{\grav}}-1\right)
 \left(\pdiffr{a}{p}\pdiffr{b}{x}\right)
 +\O\left(\hbar^2\right)\\
a\bullet b - b\bullet a &= \frac{\hbar}{\imath}
 \left(e^{-\frac{x}{\grav}}-1\right)
 \left(\pdiffr{a}{x}\pdiffr{b}{p}
 -\pdiffr{b}{x}\pdiffr{a}{p}\right)
 +\O\left(\hbar^2\right).
\end{align*}
\end{proof}

For the general quantum group formulation, we work over a general
field $k$. Since $L,V\in
\catlc{H}$ we can take their arbitrary tensor powers to define
tensor fields of arbitrary mixed rank using the same correspondence
with bicovariant bimodules. Thus $\Omega^{-1}=L\tens H$ and
$\Omega^{-1}
\tens_H\Omega^{-1}=L\tens L\tens H$ etc. We have a super-Hopf algebra
$T_{-1}(\Omega^{-1})$ and a theory of twisting of of quantum vector
fields using the same theory of Section~2. Also, since morphisms in
$\catlc{H}$ induce morphisms between bicovariant bimodules, the
evaluation map $\<\ ,\ \>:L\tens V\to k$ induces the pairing
between vector fields and 1-forms. Thus
\eqn{poisson1}{\Omega^{-1}\tens_H\Omega^1\to H,\quad \<x\tens h,v\tens g\>
=\<x,h\o\la v\>h\t g,}
\eqn{poisson2}{\Omega^{-1}\tens_H\Omega^{-1}\tens_H\Omega^1
 \to \Omega^{-1},\quad
\<x\tens y\tens h,v\tens g\>=x\tens \<y,h\o\la v\>h\t g,}
\eqn{poisson3}{\Omega^{-1}\tens_H\Omega^{-1}\tens_H
 \Omega^1\tens_H\Omega^1,\quad
\<x\tens y\tens h,v\tens w\tens g\>=\<y,h\o\la v\>\<x,h\t\la w\>h\th g,}
etc. The pairing $L\tens L\tens V\tens V\to k$ in (\ref{poisson3})
is the natural one in a braided category, namely to evaluate the
inner $L\tens V$ first and then the outer. The resulting pairing is
also the same as applying (\ref{poisson2}) to the first factor of
$\Omega^1\tens_H\Omega^1$ to obtain an element of $\Omega^{-1}
\tens_H\Omega^1$ and then applying (\ref{poisson1}).

This is not the only way to formulate vector fields (for example a more
left-right symmetric way is to consider $L\in\catrc{H}$ and $\Omega^{-1}
=H\tens L$, extending the pairing by $\<h\tens x,v\tens g\>=h\<x,v\>g$)
but it is the one natural in the context of the Woronowicz exterior
algebra (which can be viewed as based on a fixed identification of
bicovariant bimodules with $\catlc{H}$ (say)). Taking now $\Omega^n$
defined by quotients of $V^{\tens n}$ in the exterior algebra in
this approach, the natural definition of antisymmetric vector
fields is as corresponding to the appropriate subspace of $L^{\tens
n}$ dual to this quotient. In particular, the Poisson bivector
field should be an element
\[ \Pi\in \Omega^{-2}=\{x\tens y-\Psi_{L,L}(x\tens y)|\ x,y\in L\}\tens H\]
since $V\tens V$ is quotiented by $\ker(\id-\Psi_{V,V})$ in degree 2.

In general, we also need to impose a `Jacobi identity' on $\Pi$,
which can be done as follows at least in the nice case where the
quantum Poisson bracket is non-degenerate: we can consider $\Pi$ by
the above as a map $\Omega^1\to \Omega^{-1}$ and demand that it is
invertible, and that the inverse corresponds to evaluation against
some $\omega\in \Omega^{2}$ which we can demand to be closed.
Alternatively, one may attempt to develop a theory of `quantum-Lie
algebras' and use the `quantum-Lie bracket' on $L$, thereby
avoiding the invertibility assumption. This will not be attempted
here, however; for our present purposes we note that in 2
dimensions with the classical differential calculus the Jacobi
identity is redundant (similarly, every 2-form is closed). For our
particular exterior algebras the dimensions are the classical ones
(so that every 2-form is closed) and one may similarly consider any
antisymmetric bivector field as a Poisson structure. Keeping the
general form of the above classical Poisson bracket in mind, we now
give the explicit form of the quantum Poisson bracket for a natural
class of bivector fields in our example.

\begin{propos}
For the Planck scale Hopf algebra with the standard
quantum differential calculus as above, we consider $\Pi$ of the form
\[
 \Pi=(\eta^*\tens\xi^*-\Psi_{L,L}(\eta^*\tens\xi^*))\tens\pi(g)
\]
for an arbitrary function $\pi(g)$. Then the corresponding quantum
Poisson bracket is
\[
 \{a,b\}=\pi(g)\left(a_\xi b_\eta-a_\eta b_\xi+\imath A(a_\eta b_\eta
 +(a_\xi)_\eta b_\eta-(a_\eta)_\xi b_\eta)+
 (\imath A)^2 (a_\eta)_\eta b_\eta\right)
\]
where $a_\xi=\partial_\xi a$, etc. In particular,
$\pi(g)=\frac{1}{\grav}(g^{-1}-1)$ gives a deformation of the
classical Poisson structure in Proposition~\ref{prop:cbp_poisson}.
\end{propos}
\begin{proof}
We first of all use $\extd a=\xi\del_\xi a+\eta\del_\eta a$ and the
relations of the exterior algebra to obtain
\[ \extd a\wedge\extd b=\xi\wedge\eta f,
\quad f=a_\xi b_\eta-a_\eta b_\xi
 +\imath A(a_\eta b_\eta+(a_\xi)_\eta b_\eta-(a_\eta)_\xi b_\eta)+
 (\imath A)^2 (a_\eta)_\eta b_\eta.\]
Now the pairing can be computed as
\begin{align*}
\{a,b\} &=\<\Pi,\extd a\wedge \extd b\>
=\<(\eta^*\tens\xi^*-\Psi(\eta^*\tens\xi^*))\pi(g),\xi\wedge\eta f\>\\
&=\<\eta^*\tens\xi^*-\Psi(\eta^*\tens\xi^*),\xi\tens \pi(g)\o\la\eta\>
 \pi(g)\t f\\
&=\<\eta^*\tens\xi^*,\xi\tens\pi(g)\o\la\eta-\Psi_{V,V}(\xi\tens
 \pi(g)\o\la\eta)\>\pi(g)\t f\\
&=\<\eta^*\tens\xi^*,\xi\tens\pi(g)\o\la\eta-\pi(g)\o\la\eta\tens\xi\>
 \pi(g)\t f
=\<\eta^*\tens\xi^*,\xi\tens\eta\>\pi(g)f,
\end{align*}
where we used functoriality of the braiding under the evaluation
morphism to deduce
\[ \<\Psi_{L,L}(\eta^*\tens\xi^*),v\tens w)\>
 =\<\eta^*\tens\xi^*,\Psi_{V,V}(v\tens w)\>\]
for any $v,w\in V$, and then $\Psi_{V,V}(\xi\tens
w)=w\lo\la\xi\tens w
\lt=w\tens\xi$ since $\xi$ is an invariant element of the crossed module.
In the last line we used $g\la\eta=\eta+\imath A\xi$ to see that,
although $\eta$
is not invariant, the evaluation $\<Y,g^n\la\eta\>=\<Y,\eta\>$ behaves as
if it is.
In terms of functions $a(g,p)$, $b(g,p)$ we obtain
\begin{align*}
&\{:a(g,p):,:b(g,p):\}\\ &= \pi(g)
 :\left(g(g-2)a(g,p)+g(\pdiff{g}-2g+3)a(g,p+\imath A)
  +g(g-1)a(g,p+2\imath A)\right):\\
 &\qquad \bullet :g \frac{b(g,p+\imath A)-b(g,p)}{\imath A}:\\
 &\quad - \pi(g) :g\frac{a(g,p+\imath A)-a(g,p)}{\imath A}: \bullet
 :\left(g\pdiff{g} b(g,p+\imath A)+g(b(g,p+\imath A)-b(g,p))\right):.
\end{align*}
The classical limit $A\to 0$ is
\[
\{a(g,p),b(g,p)\} = \pi(g) g^2\left(\pdiffr{a}{g}\pdiffr{b}{p}
 - \pdiffr{a}{p}\pdiffr{b}{g}\right).
\]
Thus, to get the correct Poisson structure, we need
$\pi(g)=\frac{1}{\grav}(g^{-1}-1)$ (note that
$-\grav\pdiff{x}= g\pdiff{g}$).
\end{proof}

Also, if $:h:\in \pla$ is a choice of Hamiltonian then
\begin{align}
\dot{x} &=\{x,:h:\} = \frac{\grav}{\imath\hbar}:(e^{-\frac{x}{\grav}}-1)
  (h(x,p+\frac{\imath\hbar}{\grav})-h(x,p)):,\\
\dot{p} &=\{p,:h:\} = :-(e^{-\frac{x}{\grav}}-1)\pdiff{x}
  h(x,p+\imath A)):
\end{align}
are the corresponding quantum Hamilton equations of motion. For a
simple concrete example, choosing the Hamiltonian
$h(x,p)=\frac{p^2}{2m}+V(x)$ for a free particle of mass $m$ in a
potential $V(x)$, we obtain
\eqn{motion}{
 \dot{x} =\frac{1}{2 m} (e^{-\frac{x}{\grav}}-1)
  (2p-\frac{\imath\hbar}{\grav}),\quad
 \dot{p} =(e^{-\frac{x}{\grav}}-1)\pdiff{x} V(x).}
Standard quantum mechanics (i.e.\ using the commutator with $h$)
leads by contrast to
\[
\dot{x}=\frac{\imath}{\hbar}[x,h]=\frac{1}{2 m} (e^{-\frac{x}{\grav}}-1)
(2p-\frac{\imath\hbar}{\grav}e^{-\frac{x}{\grav}}),\quad
\dot{p}=\frac{\imath}{\hbar}[p,h]=(e^{-\frac{x}{\grav}}-1)\pdiff{x} V(x).
\]
Thus the quantum Hamiltonian equations of motion reduce to the
classical ones when $\hbar\to 0$ as they should, but also
approximate to the conventional quantum mechanical equations of
motion in the Planckian strongly gravitational region where
$x<<\grav$. (We recall that the quantum mechanical evolution in
this model approximates flat space when $x>>\grav$.) On the other
hand, the quantum Hamiltonian equations retain a full (quantum)
geometrical interpretation which is lost in conventional quantum
mechanics. This suggests a geometrical modification of conventional
quantum mechanics.

\section{Fourier theory on the Planck scale Hopf algebra}

In this concluding section we make some remarks about the
noncommutative Fourier theory which is known to exist on any Hopf
algebra equipped with a suitable translation-invariant integral and
a suitable exponential element. We recall first the general
formulation, which works basically when the Hopf algebra $H$ is
finite-dimensional, and in conventions suitable for our particular
example. Thus, we require $\int:H\to k$ such that $(\int h\o)h\t
=(\int h)1$ for all $h\in H$ (a right-integral) and $\int^*:H^*\to k$ such that
$\phi\o\int^*\phi\t= 1\int\phi$ for all $\phi\in H^*$ (a
left-integral), and we let $\exp=\sum e_a\tens f^a\in H\tens H^*$
denote the canonical coevaluation element (here $\{e_a\}$ is a
basis of $H$ and $\{f^a\}$ a dual basis). Then the Fourier
transform in these conventions is
\eqn{FT1}{
 \CF(h)=\left(\int e_a h\right)f^a,\quad \CF^*(\phi)=e_a \int^* f^a\phi}
and obeys
\eqn{FT2}{
 \CF\CF^*(\phi)=\antip^{-1}\phi\int e_a \int^* f^a,\quad
 \CF(h\o\<\phi,h\t\>)=\CF(h)\antip^{-1}\phi,\quad
 \CF^*(\<\phi\o,h\>\phi\t)=\antip h \CF^*(\phi).}
These elementary facts are easily proven once
one notes that $(\int g\o h)g\t=(\int g h\o)Sh\t$ for all $h,g\in H$ and a
similar identity on $H^*$. See also \cite{Ma:book},\cite{KemMa:alg} for more
discussion (and the extension to braided groups).

In our case the Planck scale Hopf algebra $\pla$ is not finite-dimensional
and there is, moreover, no purely algebraic integral. For a full treatment
one needs to introduce a Hopf-von Neumann algebra setting along the lines in
\cite{Ma:hop} and work with the integral as a weight, or one has to work with
a $C^*$ algebra setting extended to include unbounded operators. Both of these
are nontrivial and beyond our scope here. However, the bicrossproduct
form of the Hopf
algebra allows one to identify elements as normal ordered versions
of ordinary functions $f(x,p)$ and thereby to reduce integration to
ordinary integration of ordinary functions, for any class of
functions and any topological setting to which the normal ordering
extends. Therefore in this section we will initially work formally
with $x,p$ as generators (unlike the algebraic setting in the
preceeding sections) and proceed to consider formal power series in
them; however, what we arrive at in this way is a well-defined
deformed Fourier theory on functions on $\R^2$ of suitably rapid
decay, {\em motivated} by the Hopf algebra $\pla$ and consistent
with any operator algebra setting to which normal ordering extends.
This is what we shall outline in this section.

First of all, the bicrossproduct form of the Hopf algebra implies that
\[ \int :f(x,p):=\int_{-\infty}^\infty\int_{-\infty}^\infty
 \extd x\,\extd p\,f(x,p)
\]
is a {\em left}-integral on $\pla$. This is also evident
from the explicit form of the right-invariant derivatives 
(\ref{rightderiv1})-(\ref{rightderiv2}),
from which we see that the
integrals of $\del_x
:f:$ and $\del_\eta:f:$ vanish for suitably decaying $f$. On the
other hand the right-integral desired in our preferred conventions
for the Fourier theory can be similarly obtained using the
left-invariant partial differentials
(\ref{leftdel1})-({\ref{leftdel2}) stated at the end of Section~4.1
one finds
\eqn{int}{ \int :f(x,p):= \int_{-\infty}^\infty\int_{-\infty}^\infty
\extd x\,\extd p\, e^{\frac{x}{\grav}} f(x,p),}
which is the right-integral that we shall use. (Although apparently
more complicated, the resulting Fourier theory turns out to be more
computable in these conventions.)

Next, we recall from \cite{Ma:pla}\cite{Ma:book} that the
Planck scale Hopf algebra is essentially self-dual. More precisely,
if we let $\bar x,\bar p$ be dual to the $p,x$ generators in the
sense $\<\bar x,x^n p^m\>=\imath\,\delta_{n,0}\delta_{m,1}$ and $\<\bar p,
x^n p^m\>=\imath\,\delta_{n,1}\delta_{m,0}$, we have an algebraic model of
the dual of $\pla$ as
\[ \C[\bar p]\cobicross_{\frac{1}{\hbar},\frac{\grav}{\hbar}} \C[\bar x]
\subseteq (\pla)^*,\]
where
\[ [\bar p,\bar x]=\frac{\imath}{\hbar}(1-e^{-\bar x\frac{\hbar}{\grav}}),
\quad \Delta \bar x=\bar x\tens 1+1\tens\bar x,\quad \Delta\bar p=\bar p\tens 1
+e^{-\bar x\frac{\hbar}{\grav}}\tens\bar p.\] This is has the same
form as $\pla$ but with different parameter values and with the
opposite product and opposite coproduct. On this Hopf algebra we
define normal ordering as putting all the $\bar x$ to the right and
the corresponding left-integral is
\eqn{int*}{
 \int^* :f(\bar p,\bar x):=\int_{-\infty}^\infty\int_{-\infty}^\infty
\diff\bar p\,\diff\bar x\, e^{\bar x\frac{\hbar}{\grav}}
  f(\bar p,\bar x).}
Also from the
bicrossproduct form, the canonical element is \cite{Ma:book}
\eqn{exp}{\exp= \sum_{n,m}\frac{1}{n! m!\, \imath^{n+m}}
 x^np^m\tens {\bar p}^n {\bar x}^m.}
Finally, we will need explicitly the actions \cite{Ma:book}
\eqn{actions}{p\la f(x)=\imath\hbar(e^{-\frac{x}{\grav}}-1)
 \frac{\del}{\del x} f,\quad
 f(\bar x)\ra\bar p
=\frac{\imath}{\hbar}(e^{-\frac{\hbar}{\grav}\bar{x}}-1)\pdiff{\bar{x}} f}
in the bicrossproduct construction and its dual.

\begin{propos}
The quantum Fourier transform on the Planck scale Hopf
algebra is
\[
\CF(:f(x,p):)= \int_{-\infty}^\infty\int_{-\infty}^\infty \diff x\,\diff
p\, e^{-\imath(\bar{p}+\frac{\imath}{\grav})\cdot x}
 e^{-\imath \bar x\cdot (p+p\act)} f(x,p)
\]
and its dual is
\[
\CF^*(:f(\bar p,\bar x):)=\int_{-\infty}^\infty\int_{-\infty}^\infty\diff
\bar p\,\diff\bar x\, e^{-\imath\bar p\cdot x} e^{-\imath\bar x\cdot p}
 f(\ract\bar p+\bar p,\bar{x}) e^{\frac{\hbar}{\grav}\bar x},
\]
where $p\la$ acts only on the functions in $x$ to the right in the
integral ($\ract\bar p$ acts only on functions in $\bar x$ to the left).
\end{propos}
\begin{proof}
We use the reordering equality
\[ :f(p): :h(x):= :e^{p\la\cdot \pdiff{p}} h(x) f(p):
 = :f(p+p\act) h(x):\]
in $\pla$ for functions $f,h$ ($p\la$ only acts on functions of $x$). This
follows from the relation
$[p,f(x)]=p\act f(x)$ for functions $f(x)$, which is the semidirect
product form of the algebra in the bicrossproduct.
Hence,
\begin{align*}
\CF(:f(x,p):)&= \sum_{n,m}\frac{1}{n! m!\, \imath^{n+m}}
  {\bar p}^n {\bar x}^m
 \int x^n p^m :f(x,p):\\
 &=\sum_{n,m}\frac{1}{n! m!\, \imath^{n+m}} {\bar p}^n {\bar x}^m\int
   :x^n (p+p\act)^m f(x,p):\\
 &= \sum_{n,m}\frac{1}{n! m!\, \imath^{n+m}}
  {\bar p}^n {\bar x}^m\int \diff x\,\diff p\,
  e^{\frac{x}{\grav}} x^n (p+p\act)^m f(x,p)\\
 &= \int \diff x\,\diff p\, e^{-\imath(\bar p+\frac{\imath}{\grav})\cdot x}
  e^{-\imath\bar x\cdot (p+p\act)} f(x,p),
\end{align*}
where $p\la $ only acts in the powers of $x$ to its right.
In $\copla$ we have similarly
\[:f(\bar x): :h(\bar p):
 = :f(\bar x) e^{\ract\bar p\cdot \pdiff{\bar p}} h(\bar p):
 = :f(\bar x) h(\ract\bar p+p):.\]
Hence,
\begin{align*}
\CF^*(:f(\bar p,\bar x):)&=\sum_{n,m}\frac{1}{n! m!\,\imath^{n+m}}x^n p^m
 \int^* {\bar p}^n {\bar x}^m :f(\bar p,\bar x):\\
&=\sum_{n,m}\frac{1}{n! m!\, \imath^{n+m}}x^n p^m
 \int^* :\bar{p}^n \bar{x}^m f(\ract\bar p+\bar p,\bar x):\\
&=\sum_{n,m}\frac{1}{n! m!\, \imath^{n+m}}x^n p^m
 \int \diff\bar p\,\diff\bar x\, \bar{p}^n \bar{x}^m
 f(\ract\bar p+\bar p,\bar x) e^{\frac{\hbar}{\grav}\bar x}\\
&=\int \diff\bar p\,\diff\bar x\, e^{-\imath\bar{p}\cdot x}
 e^{-\imath\bar x\cdot p} f(\ract\bar{p}+\bar p,\bar x)
 e^{\frac{\hbar}{\grav}\bar x}.
\end{align*}
\end{proof}

From the properties of the Fourier transform, we see in particular
that it turns the (left-invariant) derivatives $\bar\partial_x$ and
$\bar\partial_\eta$  in (\ref{leftdel1})--(\ref{leftdel2}) into
multiplication by the corresponding element of the dual. Also,
these derivatives become right-handed derivatives $\del_{\bar x}$
and $\del_{\bar\eta}$ on $\copla$ by identifying it with the
opposite algebra and coalgebra to
$\C[x]\cobicross_{\frac{1}{\hbar},\frac{\grav}{\hbar}}\C[p]$ and
making the corresponding notational and parameter changes.

\begin{propos}
\label{prop:ft_intertwine}
\begin{gather*}
 \CF(\bar\partial _x a)=\CF(a) \imath\bar{p} e^{\frac{\hbar}{\grav}\bar{x}}
 \qquad
 \CF(\bar\partial _\eta a)=\CF(a) \frac{\imath\grav}{\hbar}
 (e^{\frac{\hbar}{\grav}\bar{x}}-1)
\\
 \CF^*(\partial_{\bar{x}}\phi)=\imath p e^{\frac{x}{\grav}}\CF^*(\phi)
 \qquad
 \CF^*(\partial_{\bar{\eta}}\phi)=\imath\grav(e^{\frac{x}{\grav}}-1)
 \CF^*(\phi)
\\
 \CF\CF^*(\phi)=(2\pi)^2\antip^{-1}\phi
\end{gather*}
\end{propos}
\begin{proof} This is a short computation to identify the partial derivatives
as $\bar\partial_x(a)=a\o\<-\imath\bar{p},a\t\>$ and
$\bar\partial_\eta(a)=a\o\<\frac{\imath\grav}{\hbar}
(e^{-\frac{\hbar}{\grav}\bar{x}}-1),a\t\>$, i.e. to identify the
corresponding elements of $L$. Similarly, $\partial_{\bar{x}}$
corresponds to $-\imath p$ and $\partial_{\bar{\eta}}$ corresponds to
$\imath\grav(e^{-\frac{x}{\grav}}-1)$ via the right coregular
action. One can then verify the analogue of (\ref{FT2}) directly in
our setting for functions of suitably rapid decay.
\end{proof}

Note that when we take the limit $\hbar\to 0$ the Hopf algebra
$\copla$ becomes $U(b_+)$ or $\kappa$-Minkowski space \cite{MaRue:bic}
with the relations
\[
 [\bar p,\bar x]=\frac{\imath}{\grav}\bar x
\]
(i.e. $\kappa=\frac{\grav}{\imath}$) regarded as a noncommutative
space. Thus,
\begin{corol}
In the classical limit $\hbar\to 0$ the Fourier transform becomes
\begin{align*}
\CF:\cbp\to\ubp, &\qquad
\CF(:f(x,p):)= \int_{-\infty}^\infty\int_{-\infty}^\infty
 \diff x\,\diff p\,
 e^{-\imath(\bar{p}+\frac{1}{\kappa})\cdot x}
 e^{-\imath\bar x\cdot p} f(x,p),\\
 \CF^*:\ubp\to\cbp, &\qquad
 \CF^*(:f(\bar p,\bar x):)=
\int_{-\infty}^\infty\int_{-\infty}^\infty\diff
\bar p\,\diff\bar x\, e^{-\imath\bar p\cdot x} e^{-\imath\bar x\cdot p}
 f(\ract\bar{p}+\bar{p},\bar{x})
\end{align*}
with $f(\bar{x})\ract\bar{p}=-\frac{\bar{x}}{\kappa}\pdiff{\bar{x}}
f$. Moreover,
\begin{gather*}
 \CF(\bar\partial_x a)=\CF(a) \imath\bar{p},
 \qquad
 \CF(\bar\partial_\eta a)=\CF(a) \imath\bar{x},\\
  \bar\partial_x (:f(x,p):)=
  :\pdiff{x} f(x,p) + \frac{\imath p}{\kappa}\,\pdiff{p} f(x,p):, \qquad
 \bar\partial_\eta(:f(x,p):)=:\pdiff{p} f(x,p):.\\
\end{gather*}
\end{corol}

The intertwiner properties of $\CF^*$ in this limit are read from
Proposition~\ref{prop:ft_intertwine} while required right-derivatives simplify to
\eqn{kappadel1}{
\partial_{\bar{x}}(:f(\bar{p},\bar{x}):)=
  :\pdiff{\bar{x}} f(\bar{p},\bar{x}):,\qquad
 \partial_{\bar{\eta}}(:f(\bar{p},\bar{x}):)=
  :-\kappa (f(\bar{p}-\frac{1}{\kappa},\bar{x})-f(\bar{p},\bar{x})):.}

We also have a dual limit $\hbar,\grav\to\infty$ with
$\frac{\grav}{\imath\hbar}=\kappa$ constant, where $\pla$ becomes
$U(\cb_-)$ (the opposite Lie algebra to $\cb_+$) and $\copla$
becomes $C(B_-)$. We regard the former as another version of
$\kappa$-Minkowski space (with opposite commutation relations).

\begin{corol}
In the limit $\hbar,\grav\to\infty$ with
$\frac{\grav}{\imath\hbar}=\kappa$ the Fourier transform becomes
\begin{align*}
\CF:U(\cb_-)\to C(B_-), &\qquad
\CF(:f(x,p):)= \int_{-\infty}^\infty\int_{-\infty}^\infty
 \diff x\,\diff p\,
 e^{-\imath\bar{p}\cdot x} e^{-\imath\bar x\cdot (p+p\act)} f(x,p),\\
 \CF^*:C(B_-)\to U(\cb_-), &\qquad
 \CF^*(:f(\bar p,\bar x):)=
\int_{-\infty}^\infty\int_{-\infty}^\infty\diff
\bar p\,\diff\bar x\, e^{-\imath\bar p\cdot (x+\frac{1}{\kappa})}
 e^{-\imath\bar x\cdot p}  f(\bar{p},\bar{x})
\end{align*}
with $p\act f(x)=-\frac{1}{\kappa}\pdiff{x} f$. Moreover,
\begin{gather*}
 \CF^*(\partial_{\bar{x}} \phi)=\imath p \CF^*(\phi),
 \qquad
 \CF^*(\partial_{\bar{\eta}} \phi)=\imath x \CF^*(\phi),\\
  \partial_{\bar{x}}(:f(\bar{x},\bar{p}):)=
  :\pdiff{\bar{x}} f(\bar x,\bar p)
   + \frac{\imath p}{\kappa}\,\pdiff{\bar{p}} f(\bar x,\bar p):, \qquad
 \partial_{\bar{\eta}}(:f(\bar x,\bar p):)=:\pdiff{\bar{p}}
  f(\bar x,\bar p):.
\end{gather*}
\end{corol}

In this limit the intertwiner properties of $\CF$ do not simplify
(we refer to Proposition~\ref{prop:ft_intertwine}), but the
corresponding derivatives become
\eqn{kappadel2}{
 \bar\partial_x (:f(x,p):)= :\pdiff{x} f(x,p):, \qquad
 \bar\partial_\eta(:f(x,p):)=:-\kappa(f(x,p-\frac{1}{\kappa})
  -f(x,p)):.}

Therefore we obtain in fact two versions of Fourier theory on
$\kappa$-Minkowski space as two limits of Fourier theory on the
Planck scale Hopf algebra. This Hopf algebra, being of self-dual
form, has the power to become both a classical but curved phase
space (the classical limit) and its dual (the second limit), in
addition to the flat space quantum mechanics limit.

There are many further possible developments of the geometry and
Fourier theory on the noncommutative phase space in this toy model
of Planck scale physics, among them quantum field theory (second
quantisation) in a first-order formalism. There is also a physical 
interpretation of the self-duality as an observable-state 
duality \cite{Ma:pla}\cite{Ma:ista}
which should be related to the noncommutative geometric picture above. 
Finally, we note that there are higher dimensional models of the 
bicrossproduct form \cite{Ma:mat}\cite{Ma:hop} which could be 
investigated from a similar point of view. These are some directions
for further work.

\section*{Acknowledgements}

R.O.\ thanks the German Academic Exchange Service (DAAD) and EPSRC for
financial support.

\appendix
\section{Direct proofs for crossed modules}

Theorem~\ref{thm:ctwist_lccat} and Proposition~\ref{prop:alpha} were deduced
somewhat indirectly
from our twisting results on bicovariant bimodules and exterior
algebras. On the other hand crossed modules $\catlc{H}$ have been
used in a variety of other contexts not related to differential
calculi and full direct proofs using conventional Hopf algebra
methods may also be useful. For completeness, we provide these
here.

\bigskip\noindent{\bf Proof of Theorem~\ref{thm:ctwist_lccat}} First, we show that
$\cfc$ is a functor, then we verify that it is monoidal. We proceed to
check the braiding and finally show that $\cfc$ is an isomorphism.

(a) $\beta_\chi$ is a coaction of $H_\chi$:
\begin{align*}
(\cou\tensor\id)\circ\beta_\chi(v)
 &=\cou(\chi\uo (\chi\umo \act v)\lo \chi\umt )
  \chi\ut \act (\chi\umo \act v)\lt  \\
 &=\cou((\chi\umo \act v)\lo \chi\umt )
  (\chi\umo \act v)\lt  \\
 &=\cou(v\lo ) v\lt =v
\end{align*}
We used the counitality of $\chi$ and $\chi^{-1}$.
\begin{align*}
&(\cop_\chi\tensor\id)\circ\beta_\chi(v)\\
 &=\cop_\chi(\chi\uo (\chi\umo \act v)\lo \chi\umt )
 \tensor\chi\ut \act (\chi\umo \act v)\lt \\
 &=\tchi\uo \chi\uo {}\o (\chi\umo \act v)\lo {}\o
  \chi\umt {}\o \tchi\umo \\
 &\quad \tensor\tchi\ut \chi\uo {}\t
  (\chi\umo \act v)\lo {}\t \chi\umt {}\t \tchi\umt
  \tensor\chi\ut \act (\chi\umo \act v)\lt \\
 &=\tchi\uo \chi\uo {}\o (\chi\umo \act v)\lo
  \chi\umt {}\o \tchi\umo \\
 &\quad \tensor\tchi\ut \chi\uo {}\t
  (\chi\umo \act v)\lt\lo
 \chi\umt {}\t \tchi\umt
  \tensor\chi\ut \act (\chi\umo \act v)\lt\lt\\
 &=\tchi\uo \chi\uo {}\o
   (\chi\umo {}\o \tchi\umo \act v)\lo
   \chi\umo {}\t \tchi\umt \\
 &\quad \tensor\tchi\ut \chi\uo {}\t
   (\chi\umo {}\o \tchi\umo \act v)\lt\lo
   \chi\umt
  \tensor\chi\ut \act
   (\chi\umo {}\o \tchi\umo \act v)\lt\lt\\
 &=\tchi\uo \chi\uo {}\o
   \chi\umo {}\o (\tchi\umo \act v)\lo
   \tchi\umt \\
 &\quad\tensor \tchi\ut \chi\uo {}\t
   (\chi\umo {}\t \act(\tchi\umo \act v)\lt )\lo
   \chi\umt
  \tensor \chi\ut \act
   (\chi\umo {}\t \act(\tchi\umo \act v)\lt )\lt \\
 &=\chi\uo
   \chi\umo {}\o (\tchi\umo \act v)\lo
   \tchi\umt \\
 &\quad\tensor \tchi\uo \chi\ut {}\o
   (\chi\umo {}\t \act(\tchi\umo \act v)\lt )\lo
   \chi\umt
  \tensor \tchi\ut \chi\ut {}\t \act
   (\chi\umo {}\t \act(\tchi\umo \act v)\lt )\lt \\
 &=\chi\uo
   \chi\umo {}\o (\tchi\umo \act v)\lo
   \tchi\umt \\
 &\quad\tensor \tchi\uo (\chi\ut {}\o
   \chi\umo {}\t \act(\tchi\umo \act v)\lt )\lo
   \chi\ut {}\t \chi\umt
  \tensor \tchi\ut \act(\chi\ut {}\t
   \chi\umo {}\t \act(\tchi\umo \act v)\lt )\lt \\
 &=\chi\uo (\tchi\umo \act v)\lo \tchi\umt \\
 &\quad\tensor \tchi\uo (\chi\umo \chi\ut \act
  (\tchi\umo \act v)\lt )\lo \chi\umt
 \tensor\tchi\ut \act (\chi\umo \chi\ut \act
  (\tchi\umo \act v)\lt )\lt \\
 &=\chi\uo (\tchi\umo \act v)\lo \tchi\umt
 \tensor\tchi\uo (\chi\umo \act \chi\ut \act
  (\tchi\umo \act v)\lt )\lo \chi\umt \\
 &\quad\tensor\tchi\ut \act (\chi\umo \act \chi\ut \act
  (\tchi\umo \act v)\lt )\lt \\
 &=\chi\uo (\tchi\umo \act v)\lo \tchi\umt
 \tensor\beta_\chi(\chi\ut \act (\tchi\umo \act v)\lt )\\
 &=(\id\tensor\beta_\chi)\circ\beta_\chi(v)
\end{align*}
We used the crossed module property (\ref{crossmod}) and the cocycle
identity (\ref{chi}). $\tchi$ denotes a second copy of $\chi$.

(b) $\beta_\chi$ together with the action obeys the crossed module
property in the twisted category:
\begin{align*}
&h\ao v\alo\tensor h\at\act v\alt\\
 &=h\ao\chi\uo (\chi\umo \act v)\lo \chi\umt
  \tensor h\at \chi\ut \act (\chi\umo \act v)\lt \\
 &=\tchi\uo h\o \tchi\umo
   \chi\uo (\chi\umo \act v)\lo \chi\umt
  \tensor \tchi\ut h\t \tchi\umt
   \chi\ut \act (\chi\umo \act v)\lt \\
 &=\tchi\uo h\o
   (\chi\umo \act v)\lo \chi\umt
  \tensor \tchi\ut h\t
   \act (\chi\umo \act v)\lt \\
 &=\tchi\uo  (h\o \act (\chi\umo \act v))\lo  h\t \chi\umt
  \tensor \tchi\ut \act(h\o \act (\chi\umo \act v))\lt \\
 &=\tchi\uo  (h\o \chi\umo \act v)\lo  h\t \chi\umt
  \tensor \tchi\ut \act(h\o \chi\umo \act v)\lt \\
 &=\tchi\uo
   (\tchi\umo \chi\uo h\o \chi\umo \act v)\lo
   \tchi\umt \chi\ut h\t \chi\umt
  \tensor \tchi\ut \act
   (\tchi\umo \chi\uo h\o \chi\umo \act v)\lt \\
 &=\tchi\uo
   (\tchi\umo \act(h\ao\act v))\lo
   \tchi\umt h\at
  \tensor \tchi\ut \act
   (\tchi\umo \act(h\ao\act v))\lt \\
&= (h\ao\act v)\alo h\at\tensor (h\ao\act v)\alt
\end{align*}
We used the crossed module property in the untwisted category and
subscripts $\ao$ etc., for the twisted coproduct and twisted
coaction.

(c) To conclude that ${\cfunct}_\chi$ is a functor, we have
to show that it maps morphisms to morphisms. Morphisms are
module-comodule maps. It is clear that ${\cfunct}_\chi$ maps module
maps to module maps since it does not alter the action. It is also
easy to see that it maps module-comodule maps to comodule maps.
Say $f:V\to W$ is a morphism in
$\catlc{H}$. Then
\begin{align*}
\beta_\chi(f(v))&=\chi\uo (\chi\umo \act f(v))\lo \chi\umt
 \tensor\chi\ut \act (\chi\umo \act f(v))\lt \\
&=\chi\uo (f(\chi\umo \act v))\lo \chi\umt
 \tensor\chi\ut \act (f(\chi\umo \act v))\lt \\
&=\chi\uo (\chi\umo \act v)\lo \chi\umt
 \tensor\chi\ut \act f((\chi\umo \act v)\lt )\\
&=\chi\uo (\chi\umo \act v)\lo \chi\umt
 \tensor f(\chi\ut \act (\chi\umo \act v)\lt )\\
&=(\id\tensor f)\circ\beta_\chi(v),
\end{align*}
as required.

(d) We proceed to show that ${\cfunct}_\chi$ is monoidal. The associativity
property of $c_\chi$ clearly reduces to $\chi$ a cocycle, and invertibility
reduces to $\chi$ invertible. Naturality of $c_\chi$ is also immediate
from its stated form. It remains to verify
that $c_\chi:\cfunct_\chi(V)\odot_\chi\cfunct_\chi(W)\to
 \cfunct_\chi(V\odot W)$
is indeed a morphism in $\catlc{H_\chi}$. For clarity we denote the tensor
product in $\catlc{H}$ by $\odot$ and that in $\catlc{H_\chi}$ by $\odot_\chi$.
For the action of $H_\chi$ (which coincides with that of $H$) we have
\begin{align*}
h\act\cnatd(v\odot_\chi w) &=h\act(\chi\umo \act v\odot \chi\umt
\act w)\\
&=h\o \chi\umo \act v\odot h\t \chi\umt \act w\\ &=\cnatd(\chi\uo
h\o \chi\umo \act v
 \odot_\chi \chi\ut h\t \chi\umt \act w)\\
&=\cnatd(h\ao\act v\odot_\chi h\at\act w)\\
&=\cnatd(h\act (v\odot_\chi w)).
\end{align*}
For the coaction, we have
\begin{align*}
&\beta_\chi\circ\cnatd(v\odot_\chi w)\\ &=\beta_\chi(\tchi\umo \act
v\odot \tchi\umt \act w)\\ &=\chi\uo (\chi\umo \act
  (\tchi\umo \act v\odot \tchi\umt \act w))\lo \chi\umt \\
 &\quad\tensor\chi\ut \act (\chi\umo \act
  (\tchi\umo \act v\odot \tchi\umt \act w))\lt \\
&=\chi\uo (\chi\umo {}\o \tchi\umo \act v
  \odot \chi\umo {}\t \tchi\umt \act w)\lo \chi\umt \\
 &\quad\tensor\chi\ut \act (\chi\umo {}\o \tchi\umo \act v
  \odot \chi\umo {}\t \tchi\umt \act w)\lt \\
&=\chi\uo (\chi\umo {}\o \tchi\umo \act v)\lo
  (\chi\umo {}\t \tchi\umt \act w)\lo \chi\umt \\
 &\quad\tensor\chi\ut \act
  ((\chi\umo {}\o \tchi\umo \act v)\lt
  \odot (\chi\umo {}\t \tchi\umt \act w)\lt )\\
&=\chi\uo (\chi\umo {}\o \tchi\umo \act v)\lo
  (\chi\umo {}\t \tchi\umt \act w)\lo \chi\umt \\
 &\quad\tensor\chi\ut {}\o \act
  (\chi\umo {}\o \tchi\umo \act v)\lt
  \odot \chi\ut {}\t \act
  (\chi\umo {}\t \tchi\umt \act w)\lt \\
&=\chi\uo (\chi\umo \act v)\lo
  (\chi\umt {}\o \tchi\umo \act w)\lo
  \chi\umt {}\t \tchi\umt \\
 &\quad\tensor\chi\ut {}\o \act
  (\chi\umo \act v)\lt
  \odot \chi\ut {}\t \act
  (\chi\umt {}\o \tchi\umo \act w)\lt \\
&=\chi\uo (\chi\umo \act v)\lo
  \chi\umt {}\o (\tchi\umo \act w)\lo
  \tchi\umt \\
 &\quad\tensor\chi\ut {}\o \act
  (\chi\umo \act v)\lt
  \odot \chi\ut {}\t
  \chi\umt {}\t \act(\tchi\umo \act w)\lt \\
&=\chi\uo (\chi\umo \act v)\lo
  \chi\umt {}\o \hchi\umo \tchi\uo
  (\tchi\umo \act w)\lo
  \tchi\umt \\
 &\quad\tensor\chi\ut {}\o \act
  (\chi\umo \act v)\lt
  \odot \chi\ut {}\t
  \chi\umt {}\t \hchi\umt \tchi\ut
  \act(\tchi\umo \act w)\lt \\
&=\chi\uo (\chi\umo {}\o \hchi\umo \act v)\lo
  \chi\umo {}\t \hchi\umt \tchi\uo
  (\tchi\umo \act w)\lo
  \tchi\umt \\
 &\quad\tensor\chi\ut {}\o \act
  (\chi\umo {}\o \hchi\umo \act v)\lt
  \odot \chi\ut {}\t
  \chi\umt \tchi\ut
  \act(\tchi\umo \act w)\lt \\
&=\chi\uo \chi\umo {}\o (\hchi\umo \act v)\lo
  \hchi\umt \tchi\uo
  (\tchi\umo \act w)\lo
  \tchi\umt \\
 &\quad\tensor\chi\ut {}\o
  \chi\umo {}\t \act(\hchi\umo \act v)\lt
  \odot \chi\ut {}\t
  \chi\umt \tchi\ut
  \act(\tchi\umo \act w)\lt \\
&=\chi\uo (\hchi\umo \act v)\lo
  \hchi\umt \tchi\uo
  (\tchi\umo \act w)\lo
  \tchi\umt \\
 &\quad\tensor\chi\umo
  \chi\ut \act(\hchi\umo \act v)\lt
  \odot \chi\umt \tchi\ut
  \act(\tchi\umo \act w)\lt \\
&=(\id\tensor\cnatd)(\chi\uo (\hchi\umo \act v)\lo
  \hchi\umt \tchi\uo
  (\tchi\umo \act w)\lo
  \tchi\umt \\
 &\quad\tensor
  \chi\ut \act(\hchi\umo \act v)\lt
  \odot_\chi \tchi\ut
  \act(\tchi\umo \act w)\lt )\\
&=(\id\tensor\cnatd)(v\alo w\alo
 \tensor v\alt\odot_\chi w\alt)\\
&=(\id\tensor\cnatd)\circ\beta_\chi(v\odot_\chi w).
\end{align*}
We used the crossed module property and the cocycle identities as
before, and $\hchi$ denotes a third copy of $\chi$.

(e) We next show that $\cfunct_\chi$ preserves the
braiding. Thus,
\begin{align*}
\cnatd\circ\Psi_\chi(v\odot_\chi w)
 &=\cnatd(v\alo\act w\odot_\chi v\alt)\\
 &=\cnatd(\chi\uo (\chi\umo \act v)\lo \chi\umt \act w
 \odot_\chi \chi\ut \act (\chi\umo \act v)\lt )\\
 &=(\chi\umo \act v)\lo \chi\umt \act w
 \odot (\chi\umo \act v)\lt \\
 &=\Psi(\chi\umo \act v\odot \chi\umt \act w)\\
 &=\Psi\circ\cnatd(v\odot_\chi w).
\end{align*}

(f) It remains to be shown that $\cfc$ is an isomorphism. The inverse
operation to the twisting by $\chi$ is twisting by
$\chi^{-1}$; we verify that the coaction twisted by $\chi$ and
then twisted by $\chi^{-1}$ is the original coaction. Thus,
\begin{align*}
(\beta_\chi)_{\chi^{-1}}(v) &=\chi\umo (\chi\uo \act v)\alo\chi\ut
 \tensor\chi\umt \act (\chi\uo \act v)\alt\\
&=\chi\umo \tchi\uo (\tchi\umo \chi\uo \act v)\lo
  \tchi\umt \chi\ut
 \tensor\chi\umt \tchi\ut \act
  (\tchi\umo \chi\uo \act v)\lt \\
&=v\lo \tensor v\lt =\beta(v).
\end{align*}
For the monoidal structure, one sees immediately
that $\cnat_{\chi^{-1}}\circ\cnat_\chi$ is the identity transformation.

Theorem~\ref{thm:ptwist_lccat} can likewise be proven directly or
else be obtained by
dualisation of Theorem~\ref{thm:ctwist_lccat} using conventional
methods.

\bigskip
\noindent{\bf Proof of the intertwiner property of $\alpha$ in
Proposition~\ref{prop:alpha} in $H_\chi$ setting}
We give the result here in the same coproduct twist setting as
Theorem~\ref{thm:ctwist_lccat} proven above (the version in
Proposition~\ref{prop:alpha} is the
dual of this and can be obtained by the same methods or by
dualisation of the proof.) Thus, for a Hopf algebra $H$ viewed in
$\catlc{H}$ by the coproduct (the regular coaction) and adjoint
action, and a cocycle $\chi\in H\tens H$ we show that
\eqn{alpha*}{ \alpha:\cfunct_\chi(H)\to H_\chi,\quad \alpha(h)
=(\chi\umo\la h)\chi\umt}
is an isomorphism of crossed modules, where $H_\chi$ is viewed in
$\catlc{H_\chi}$ by its coproduct $\Delta_\chi$ and its adjoint
action. In fact, that the actions are intertwined is known from
\cite{GurMa:bra} (in another context) so we need only to show that
the coactions are intertwined.

On the one hand, writing $\beta_\chi$ for the coaction induced by
Theorem~\ref{thm:ctwist_lccat} on $\cfunct_\chi(H)$ by transforming the regular
coaction, we have
\[(\id\tens\alpha)\beta_\chi(h)=\chi\uo(\chi\umo\la h)\o\chi\umt
\tens(\chi'\umo\chi\ut\la(\chi\umo\la h)\t)\chi'\umt .\]
We require this to coincide with
\begin{align*}
\Delta_\chi\alpha(h)&=\chi\uo(\chi\umo\la h)\o\chi\umt\o\chi'\umo\tens\chi\ut
(\chi\umo\la h)\t\chi\umt\t\chi'\umt\\
&=\chi\uo(\chi\umo\o\chi'\umo\la h)\o\chi\umo\t\chi'\umt\tens\chi
\ut(\chi\umo\o\chi'\umo\la h)\t\chi\umt\\
 &=\chi\uo \chi\umo\o (\chi'\umo\la h)\o\chi'\umt\tens\chi\ut(
 \chi\umo\t\la(\chi'\umo\la h)\t)\chi\umt .
 \end{align*}
using the cocycle axiom for $\chi$ and then the crossed module axiom. 
Comparing these expressions and substituting the quantum group
adjoint action of $H$ for $\la$ we see that these expression
coincide in view of the identity
\[ \chi\uo\tens\chi\umo\o\chi\ut\o\tens (S\chi\ut\t)(S\chi\umo\t)\chi\umt
=\chi\uo\chi\umo\o\tens\chi\ut\chi\umo\t\tens (S\chi\umo\th)\chi\umt .\]
This is equivalent (by using the cocycle condition (\ref{chi}) on
the left hand side repeatedly) to
\[ \chi\uo\o\tens\chi\uo\t\tens (S(\chi'\umo\chi\ut))\chi'\umt=
\chi\umo\o\tens\chi\umo\t\tens (S\chi\umo\th)\chi\umt\]
which reduces to
\[
 \chi\uo\tens (S(\chi'\umo\chi\ut))\chi'\umt=
\chi\umo\tens (S\chi\umo\th)\chi\umt .
\]
This identity is readily proven from the properties of
$U^{-1}=(S\chi\umo)\chi\umt$ in \cite{Ma:book} using the cocycle
condition. The inverse of the map $\alpha$ is also readily supplied
by similar means, so it forms an isomorphism of crossed modules.

\baselineskip 18pt

\begin{thebibliography}{10}

\bibitem{Ma:cla}
S.~Majid.
\newblock Classification of bicovariant differential calculi.
\newblock {\em J. Geom. Phys.}, 25:119--140, 1998.

\bibitem{BegMa:dif}
E.~Beggs and S.~Majid.
\newblock Quasitriangular and differential structures on bicrossproduct Hopf
  algebras.
\newblock {\em Preprint}, Damtp/96-97, 1996.

\bibitem{Dri:hop}
V.G. Drinfeld.
\newblock {H}opf algebras and the quantum {Y}ang-{B}axter equations.
\newblock {\em Sov. Math. Dokl.}, 32:254--258, 1985.

\bibitem{Dri:qua}
V.G. Drinfeld.
\newblock Quasi{H}opf algebras.
\newblock {\em Leningrad Math. J.}, 1:1419--1457, 1990.

\bibitem{Ma:book}
S.~Majid.
\newblock {\em Foundations of Quantum Group Theory}.
\newblock Cambridge Univeristy Press, 1995.

\bibitem{GurMa:bra}
D.I. Gurevich and S.~Majid.
\newblock Braided groups of {H}opf algebras obtained by twisting.
\newblock {\em Pac. J. Math.}, 162:27--44, 1994.

\bibitem{Ma:pla}
S.~Majid.
\newblock {H}opf algebras for physics at the {P}lanck scale.
\newblock {\em J. Classical and Quantum Gravity}, 5:1587--1606, 1988.

\bibitem{MaRue:bic}
S.~Majid and H.~Ruegg.
\newblock Bicrossproduct structure of the {$\kappa$}-{P}oincar{\'e} group and
  non-commutative geometry.
\newblock {\em Phys. Lett. B}, 334:348--354, 1994.

\bibitem{BFFLS:def}
F.~Bayen, M.~Flato, C.~Frondsdal, A.~Lichnerowicz, and D.~Sternheimer.
\newblock Deformation theory and quantisation i,ii.
\newblock {\em Ann. Phys.}, 111:61--151, 1978.

\bibitem{Oe:cla}
R.~Oeckl.
\newblock Classification of differential calculi on $U_q(b_+)$, classical
  limits, and duality.
\newblock {\em Preprint}, Damtp-1998-86, math.QA/9807097.

\bibitem{Dri}
V.G. Drinfeld.
\newblock Quantum groups.
\newblock In A.~Gleason, editor, {\em Proceedings of the {ICM}}, pages
  798--820, Rhode Island, 1987. AMS.

\bibitem{Rad:str}
D.~Radford.
\newblock The structure of {H}opf algebras with a projection.
\newblock {\em J. Algebra}, 92:322--347, 1985.

\bibitem{Yet:rep}
D.N. Yetter.
\newblock Quantum groups and representations of monoidal categories.
\newblock {\em Math. Proc. Camb. Phil. Soc.}, 108:261--290, 1990.

\bibitem{Wor:dif}
S.L. Woronowicz.
\newblock Differential calculus on compact matrix pseudogroups (quantum
  groups).
\newblock {\em Commun. Math. Phys.}, 122:125--170, 1989.

\bibitem{Brz:rem}
T.~Brzezi\'nski.
\newblock Remarks on bicovariant differential calculi and exterior {{H}opf}
  algebras.
\newblock {\em Lett. Math. Phys.}, 27:287, 1993.

\bibitem{Ker:gra}
R.~Kerner.
\newblock $Z_3$-graded algebras and the cubic root of the supersymmetry
  transformations.
\newblock {\em J. Math. Phys.}, 33:403--411, 1992.

\bibitem{Ma:euc}
S.~Majid.
\newblock {$q$}-{E}uclidean space and quantum {W}ick rotation by twisting.
\newblock {\em J. Math. Phys.}, 35:5025--5034, 1994.

\bibitem{Ma:hop}
S.~Majid.
\newblock {H}opf-von {N}eumann algebra bicrossproducts, {K}ac algebra
  bicrossproducts, and the classical {Y}ang-{B}axter equations.
\newblock {\em J. Funct. Analysis}, 95:291--319, 1991.
\newblock From PhD Thesis, Harvard, 1988.

\bibitem{Ma:mat}
S.~Majid.
\newblock Matched pairs of {L}ie groups associated to solutions of the
  {Y}ang-{B}axter equations.
\newblock {\em Pac. J. Math.}, 141:311--332, 1990.
\newblock From PhD Thesis, Harvard, 1988.

\bibitem{Ma:ista}
S.~Majid.
\newblock Duality principle and braided geometry.
\newblock In {\em Strings and Symmetries}, volume 447 of {\em Lec. Notes in
  Phys.}, pages 125--144. Springer, 1995.

\bibitem{KemMa:alg}
A.~Kempf and S.~Majid.
\newblock Algebraic $q$-integration and {F}ourier theory on quantum and braided
  spaces.
\newblock {\em J. Math. Phys.}, 35:6802--6837, 1994.

\end{thebibliography}

\end{document}